\documentclass[12pt]{article}
\usepackage{amssymb}

\makeatletter

\@addtoreset{equation}{section}

\makeatother

\newtheorem{teo}[equation]{Theorem}
\newtheorem{lema}[equation]{Lemma}
\newtheorem{sublema}[equation]{Sublemma}
\newtheorem{coro}[equation]{Corollary}
\newtheorem{prop}[equation]{Proposition}

\newtheorem{defin}[equation]{Definition}
\newtheorem{rem}[equation]{Remark}

\renewcommand\thefigure{\thesection.\@arabic\c@figure}

\renewcommand\thetable{\thesection.\@arabic\c@table}

\def\E{{\mathbb E}}
\def\P{{\mathbb P}}
\def\Q{{\mathbb Q}}
\def\R{{\mathbb R}}
\def\Z{{\mathbb Z}}

\def\N{{\mathbb N}}

\def\PP{{\mathbb P}}

\def\N{{\mathbb N}}

\def\C{{\bf C}}

\def\G{{\bf G}}

\def\sqr{\vcenter{
         \hrule height.1mm
         \hbox{\vrule width.1mm height2.2mm\kern2.18mm\vrule width.1mm}
         \hrule height.1mm}}                  
\def\square{\ifmmode\sqr\else{$\sqr$}\fi}
\def\one{{\bf 1}\hskip-.5mm}

\def\ga{{\gamma}}
\def\th{{\theta}}

\def\proof{\noindent{\bf Proof. }}
\def\A{{\bf A}}

\def\E{{\mathbb E}}
\def\P{{\mathbb P}}

\def\reff#1{(\ref{#1})}
\def\ani{\G}
\def\jota{{\bf J}}
\def\supp{{\rm Supp}\,}
\def\III{{\rm Int}}
\def\qint{\overline{\Q}}
\def\rea{\underline C}
\def\kep{\underline K}
\def\era{\underline D}

\def\mark{{\rm Mark}\,}
\def\basis{{\rm Basis}\,}
\def\dist{{\rm dist}\,}
\def\life{{\rm Life}\,}
\def\birth{{\rm Birth}\,}

\def\tl{{\rm TL}\,}
\def\sd{{\rm SD}\,}
\def\ss{{\rm SS}\,}
\def\que{\aleph}

\def\norm#1{\left\|#1\right\|}

\begin{document}

\title{{\bf Spatial birth-and-death processes\\ in random environment}}

\author{Roberto Fern\'{a}ndez\footnote{Laboratoire de Math\'{e}matiques
Rapha\"el Salem, UMR 6085 CNRS--Universit\'{e} de Rouen,
F 76821 Mont Saint Aignan--Cedex,
FRANCE; {\tt Roberto.Fernandez@univ-rouen.fr}}, Pablo A.
Ferrari\footnote{Instituto de Matem\'{a}tica e Estat\'{\i}stica -
USP, Caixa Postal 66281, CEP 05389-970 S\~ao Paulo,
BRAZIL; {\tt pablo@ime.usp.br}} and Gustavo R. Guerberoff\footnote{Instituto de Matem\'{a}tica y
Estad\'{\i}stica Rafael Laguardia, Facultad de Ingenier\'{\i}a - Universidad de la Rep\'{u}blica,
Julio Herrera y Reissig 565, CP 11300, Montevideo, URUGUAY; {\tt gguerber@fing.edu.uy}.}}

\maketitle

\begin{abstract}

  We consider birth-and-death processes of objects (animals) defined
  in $\Z^d$ having unit death rates and random birth rates. For
  animals with uniformly bounded diameter we establish conditions on
  the rate distribution under which the following holds for almost all
  realizations of the birth rates: (i) the process is ergodic with at
  worst power-law time mixing; (ii) the unique invariant measure has
  exponential decay of (spatial) correlations; (iii) there exists a
  perfect-simulation algorithm for the invariant measure.  The results
  are obtained by first dominating the process by a backwards oriented
  percolation model, and then using a multiscale analysis due to Klein
  to establish conditions for the absence of percolation.


\end{abstract}

\tableofcontents

\section{Introduction}

A number of reasons explain the standing interest in the study of
birth-and-death processes.  First, they are important probabilistic
constructions in their own right, with proven potential to generate
and test new mathematical approaches.  Second, they offer the correct
framework to formalize and study several important statistical
applications.  Among the more recent ones, we mention point
processes (Baddeley and van Lieshout, 1995\nocite{badlie95}; Strauss,
1995\nocite{str75}; Baddeley, Kendall and van Lieshout,
1996\nocite{bkl96}) and loss networks (Kelly, 1991\nocite{kel91},
Ferrari and Garcia, 1998\nocite{fg98}).  Third, they have become
useful tools for other mathematical endeavours.  Examples of this are
Kendall's (1998\nocite{ken98}) perfect simulation scheme and its
offsprings (see M{\o}ller, 2001\nocite{mol01} for a review) and the
study of hard-core statistical mechanical systems presented in
Fern\'{a}ndez, Ferrari and Garcia (1998, 2001)\nocite{ffg98,ffg01}.

While the non-random version of these processes has a well developed
theory (see, for instance, M{\o}ller and Waagepetersen,
2004\nocite{molwaa04}), little is found in the literature for
processes with random rates.  Yet, in many situations a process with
spatially homogeneous randomness is a more realistic model than a
spatially homogeneous deterministic one.  In general, the behavior of
processes with ``strong'' disorder is expected to be qualitatively
different from that of the non-disordered counterpart.  In this paper,
however, we tackle the complementary issue.  We determine sufficient
conditions for the disorder to be ``weak'' in the sense that it leads
to processes exhibiting ergodic properties not far from those of the
non-random versions.

We set our processes on a lattice that for convenience is taken to be
$\Z^d$. The objects being born are called \emph{animals}.  They are
supported on finite subsets of $\Z^d$, which, for technical reasons,
are assumed to be uniformly bounded in size.  Each animal $\gamma$ has
an associated birth rate $w^{\jota}(\gamma)$ parametrized by a random
variable $\jota$ defined on a certain probability space.  Each
realization of $\jota$ is a random environment.  In addition, there is
an \emph{incompatibility relation} between animals which is also
assumed to be finite-range.  The birth-and-death process is defined as
follows: for each fixed environment, each animal $\gamma$
attempts to appear at rate $w^{\jota}(\gamma)$ but to succeed, it must
past a test, in general of stochastic nature, involving incompatible
animals present at the moment of the attempt.  Once born, animals
disappear at a unit rate.  Our main result, Theorem \ref{155},
establishes two conditions on the disorder [hypotheses \reff{eq:r1}
and \reff{eq:r3} below] under which the process is ergodic and
space-time mixing for almost all random environments.  As in other
models evolving in the presence of frozen-disorder (Klein,
1995\nocite{kle95}; Gielis and Maes, 1996\nocite{giemae96}), while
space mixing remains exponential, time relaxation is only proven to be
faster than any power, with a bound of the form $\exp[-m\ln^q(1+t)]$,
for some $m, q>0$.

Two conditions are imposed on the disorder.  On the one hand, in
\reff{eq:r1} the average local rate of \emph{attempted} births is
asked to have some finite logarithmic moment.  This is a relatively
mild condition which is independent of the incompatibility relation.
Its role [see Lemma \ref{varsovia2} and formulas
\reff{rr:7}--\reff{eq:rr.23}] is to ensure that there is a
not-too-small probability that locally no animal is found to be alive
throughout a sufficiently thin time slice, even in the absence of
incompatibility restrictions.  The second disorder condition
\reff{eq:r3} requires a weighted birth-rate of incompatible animals to
have a sufficiently small mean.  Weights are defined by a certain
\emph{size} function $S(\gamma)$ which in principle can be chosen in
any way leading to the validity of the condition.  This is a situation
very much in the spirit of cluster-expansion formalisms; see, for
instance, Dobrushin (1996)\nocite{dob96}.  Usually ---but not
optimally--- the size function is chosen as the number of sites in the
support of the animal.  This second disorder condition guarantees that
large fluctuations in rate values are sufficiently sparse to have a
negligible effect at sufficiently large scales.

In a subsequent corollary (Corollary \ref{coro.main}), we determine a
convenient sufficient criterium for both disorder hypotheses.  We show
that, barring very exceptional settings, there exists
an
$\widetilde\varepsilon$ (which we do not try to optimize) such that if
\begin{equation}
\label{eq:rr.00}
\E\Bigl[\sup_x \sum_{\theta \ni x} \left|H(\theta)\right|
\,w^\jota(\theta)\Bigr] \;\le\; \widetilde \varepsilon
\end{equation}
the process is almost-surely ergodic and mixing.  Here
$\left|H(\gamma)\right|$ is the cardinality of the \emph{halo}, which
is the region around $\gamma$ that determines incompatibility with
other animals [see \reff{eq:rr.50}].

Our analysis is based on a construction of birth-and-death processes
described in Fern\'{a}ndez, Ferrari and Garcia (1998\nocite{ffg98},
2001\nocite{ffg01}).  In this approach, the process is constructed by
resorting first to a \emph{free} process obtained by turning off
incompatibilities and allowing every attempted birth to succeed.  This
free process corresponds to a marked Poisson process, where the marks
are the lifetime of the animal and a random variable to be used in the
compatibility test.  It is useful to visualize each realization of
this free process as the collection of cylinders determined by the
animals alive during a certain lifespan.  These cylinders inherit the
incompatibility relation of the animlas forming their sections.  Of
particular importance for the construction of the interacting process
are those incompatible cylinders that are alive when a given cylinder
is born.  They are called \emph{ancestors} of the latter.  The
interacting process can be constructed if each cylinder has only a
finite number of generations of ancestors.  In this case, one can
descend down genealogical trees deciding, by means of succesive
incompatibility tests, which cylinders of the free process remain for
the interacting one.  Formally, the relation ``being ancestor of''
defines a (backwards) oriented percolation model in the space of
cylinder configurations.  The construction scheme succeeds if this
model does not exhibit infinite percolation.

A sufficient condition for this lack of percolation is, in turns,
obtained by resorting to a dominant multitype branching process whose
subcriticality allows the construction to work.  The subcriticality
condition is, precisely,
\begin{equation}
    \label{eq:r00}
    \sup_\gamma \, \frac{1}{S(\gamma)} \sum_{\theta
{\rm\ incomp.\, with\ } \gamma} S(\theta) \,  w(\theta) \;<\; 1\;.
\end{equation}
The left-hand side of this inequality is the mean number of branches
of the dominating branching process, and acts as a driving parameter:
Time and space rates of convergence and mixing rates can be explicitly
obtained in terms of it.

This approach can not be directly applied for random birthrates
because, except in trivial cases, condition \reff{eq:r00} is violated
with probability one.  An additional argument is needed showing that
these violations are so sparse that the process is basically driven by
the behavior in the overwhelmingly present ``good'' (regular) regions.
As often in the study of disordered systems, this additional argument
takes the form of a \emph{multiscale analysis}.  In this paper we
adapt a time-tested multiscale argument whose origin can be traced to
von Dreifus' (1987) dissertation\nocite{dre87}.  The argument was
later adapted by Campanino and Klein (1991\nocite{camkle91}), and
Campanino, Klein and Perez (1991\nocite{camkleper91}) to the study of
$(d+1)$-dimensional systems with $d$-dimensional disorder, which is
our setting here.  In the present work we follow the more general
version due to Klein (1994\nocite{kle94}) which has also been the
basis of Gielis and Maes' (1996\nocite{giemae96}) study of spin-flip
dynamics.  Our proofs are similar to those of these references, and
the informed reader will recognize many points in common, such as
conditions \reff{pf1}, \reff{pf3} and \reff{pf4} below.  But there is
a number of small but frequent adaptations needed to accomodate
``blob'' rather than edge percolation.  For the sake of
clarity and completeness we have prefered to do a self-contained
exposition of the proof, rather than a catalogue of differences with
respect to preceding references.  The paper size is comparable in both
cases.  We have, nevertheless, kept most of the notation, and of the
architecture of the proof, adopted by Klein (1994) and Gielis and Maes
(1996\nocite{giemae96}) for the benefit of readers familiar with these
papers.

The initial ingredient of the multiscale scheme is a sequence of boxes
of increasing size.  The linear sizes of the boxes are the
\emph{scales}.  Sites are then classified as \emph{regular} or {\em
  singular} at a given scale according to whether there is an
appropriate decay of the integrated connectivity function to the
boundary of a box around the site.  For the procedure to succeed there
must exist a choice of scales that noticeably improves the probability
of a site to be regular at succesive scales.  This requires the
singular regions to be sparsely distributed within a sea of regular
sites; the distribution becoming more and more diluted as the
observation scale grows.  In the present setting of space-time
(backwards) percolation, boxes are space-time parallelepipeds with a
cubic spatial section.  The space scales grow as a power law, but the
time scale needs to grow faster.  This is a consequence of the
$d$-dimensional nature of the disorder.  While regular and singular
regions alternate in spatial directions, they are frozen in the time
direction.  Their effects are limited only by the deaths and lack of
births of the concerned animals, which are processes of a more
correlated nature than the (spatial) disorder distribution.  Hence,
the connectivity function decays much more slowly in time than in
space.  The time scale must satisfy two complementary requirements.
First, it must grow fast enough to ensure a faster-than-power decay in
the percolation probability from bottom to top inside a box.  But,
second, the growth must be suficiently slow to guarantee that
time-like connections be mostly established inside a box, rather than
by creeping through a sequence of spacially consecutive boxes.
Mathematically, the first requirement is used in the proof of Theorem
\ref{behavior}, to arrive to formula \reff{tiempo}, while the second
one appears explicitly in the proofs of Sublemma \ref{sl.second} and
Lemma \ref{varsovia2}.  The compromise is achieved by chosen a
function
$T(L)$ so
the bottom-top and left-right percolation probabilities inside a box
of size $L\times T(L)$ decrease at a comparable rate, namely
exponentially with $L$.  This corresponds to a stretched-exponential
dependence [(see \reff{con1}--\reff{eq:rr.20.1}].  \medskip

In statistical mechanical terms, our results are, somehow, the
low-temperature counterpart of those of Gielis and Maes
(1996\nocite{giemae96}).  While these authors study
single-spin-flipping evolutions leading
to high-temperature invariant measures, here we focus on gases of
``deffects'' characteristic of low-temperature invariant measures.
The comparison is not fully valid, however, because we are able to
study only families of uniformly bounded deffects.

\section{Definitions, examples and results}

\subsection{Basic definitions}\label{s.basic}

\paragraph{Animals}
We consider the lattice $\Z^d$, $d\ge 1$, (or, in general, the
vertices of a graph with uniformly bounded coordination number)
endowed with some norm.  For concreteness we adopt $\norm{x}
 = \max \{|x_i|: i=1,2, \ldots , d\}$, $x \in \Z^d$.  An
animal model in $\Z^d$ is defined by a countable family $\ani$ of
objects, the \emph{animals}, for which there exists a map
\begin{equation}
  \label{eq:2}
  \begin{array}{rcl}
\ani & \longrightarrow & \mathcal{P} (\Z^d)\\
\ga & \mapsto & V(\ga)
  \end{array}
\end{equation}
such that (i) $V(\ga)$ is a finite set ---whose elements are called
the \emph{vertices} or \emph{sites} of $\ga$--- and (ii) there is a
finite number of animals
asociated to each fixed set of vertices:
$V^{-1}(V(\ga))$ is finite for every $\ga\in\ani$. Animal
configurations are elements $\underline\xi$ of $\N^\ani$.

Geometry is introduced through the map $V$ in the obvious manner: The
\emph{diameter} of an animal $\ga$ is the diameter of its set
$V(\ga)$, the \emph{distance} between two animals is the distance
between the corresponding sets of vertices, etc.  To abbreviate we
shall often denote $x\in\ga$, for $x\in\Z^d$, rather than $x\in
V(\ga)$.  For a region $\Lambda\subset\Z^d$ (of finite or infinite
cardinality), we denote $\ani_\Lambda$ the set formed by all animals
$\ga$ with $V(\ga)\subset\Lambda$ (animals \emph{in} $\Lambda$).  The
corresponding configuration will get a subscript $\Lambda$:
$\underline\xi_\Lambda\in\N^{\ani_\Lambda}$.  The omission of a
subscript
$\Lambda$ indicates $\Lambda=\Z^d$. Animals
are also characterized by a \emph{size function}.  In general terms, a
size function is any funcion $S:\ani \to [1,\infty)$ that can be used
in the convergence condition \reff{eq:r3} below.  In the examples of
next section, this function is just the number of sites in $V(\gamma)$.

In this paper we suppose that there exists a number $\ell_1$ such that
\begin{equation}
  \label{eq:3}
  {\rm diam}\, (\ga) \;\le\; \ell_1 \qquad \hbox{for all }
  \ga\in\ani\;.
\end{equation}

\paragraph{Animal interactions}
We introduce an \emph{interaction function} $M(\,\cdot\,|\,\cdot\,):
\ani\times \N^{\ani} \to [0,1]$, such that $M(\ga\vert \underline \xi)$
is the probability that an attempted birth of $\ga\in\ani$ actually
occurs when the current configuration of animals is $\underline\xi$.
If this function takes only the values 0 or 1 we refer the interaction
as \emph{deterministic}.  This happens, for instance, for conditions
such as volume- or perimeter-exclusion.  The function $M$ determines
the (binary) interaction matrix
\begin{equation}
  \label{200}
  \III(\ga,\th) \;:=\; \one\Bigl\{\sup_{\underline \xi}
\Bigl|M(\ga\vert \underline \xi) - M(\ga\vert
\underline\xi+\delta_\theta)\Bigr| \neq 0\Bigr\}
\end{equation}
where $\one\{A\}$ is the indicator function of the set $A$,
$\delta_\theta$ is the configuration in $\N^\ani$ having only $\theta$
present, and $\underline \xi+\delta_{\theta}$ is the configuration
obtained by adding the animal $\theta$ to the configuration
$\underline \xi$.  We say that $\gamma$ is \emph{incompatible} with
$\theta$, and denote $\gamma \not\sim \theta$, iff $\III(\ga,\th)=1$;
otherwise we say that $\gamma$ is compatible with
$\theta$.
A family of animals $\Gamma\subset \ani$ is a
\emph{compatible family} [\emph{herd} in Dobrushin's (1996)
nomenclature\nocite{dob96}] if its elements are pairwise compatible.
We assume that there exists a finite $\ell_2\in\R$ such that
\begin{equation}
  \label{eq:11}
{\rm dist}(\ga,\th)>\ell_2 \; \Longrightarrow \;\III(\ga,\th)=0\;.
\end{equation}

\paragraph{The random environment}
Each animal $\ga$ has an associated \emph{birth rate} $w^\jota(\ga)$
which depends on a variable $\jota$, the \emph{random environment},
belonging to a certain probability space $({\cal J}, \PP)$.  We assume
that the \emph{disorder measure} $\PP$ is such that
\begin{equation}
  \label{eq:10}
\Gamma \hbox{ compatible } \Longrightarrow \
\{w^\jota(\gamma) \,:\, \gamma\in\Gamma\} \ \hbox{independent.}
\end{equation}

\paragraph{Birth-and-death processes in random environment}
For each fixed environment $\jota$ and each $\Lambda\subset\Z^d$ we
consider the interacting birth-and-death processes
formally defined by the generator:
\begin{equation}
  \label{gen}
  A_\Lambda^\jota F(\underline\eta) \;=\; \sum_{\gamma \in \ani_\Lambda}
  \Bigl\{M(\ga|\underline\eta) \,w^{\jota}(\ga)\,
  [F(\underline\eta+\delta_\ga)- F(\underline\eta)]\;+\;  \eta(\ga)\,
  [F(\underline\eta-\delta_\ga)- F(\underline\eta)]  \Bigr\}\;,
\end{equation}
where $F$ is a real continuous function on $\N^{\ani_\Lambda}$.

In words, \reff{gen} says that when the current configuration of
animals is $\underline\eta$, each animal $\ga$ appears at rate $M(\ga
|\underline \eta)\,w^{\jota}(\ga)$ and disappears at rate 1 (if it is
present).  This is an interacting birth-and-death process of animals
with quenched disorder.  The factor $w^{\jota}(\ga)$ acts as the rate
of an internal Poissonian clock marking possible birth instants.  The
actual birth takes place only if a further test, determined by
$M(\ga|\,\cdot\,)$, is passed.  From the definition \reff{200} of the
interaction matrix we see that if $\III(\ga,\th)$ takes the value 1
(resp.\ 0) the presence of $\th$ may have (resp.\ does not have) an
influence on the birth rate of $\ga$.

Unit death-rates are no loss of generality.  In general, our
$w^\jota(\gamma)$ stands for the ratio ``birth rate / death rate'' of
$\gamma$.  In fact, a similar treatment is possible for
interacting-birth interacting-death processes.  These are processes
with a generator as in \reff{gen} but with a factor $M_{\rm
  death}(\ga|\underline\eta)\,w^\jota_{\rm death}(\gamma)$ multiplying
the last summand.  Our theory can be easily adapted if the death rates
$w^\jota_{\rm death}(\gamma)$ are bounded below by a strictly positive
number, uniformly in $\gamma$ and $\underline\eta$.

In general, there may be forbidden cases for which the matrix $M$
takes value 0.  An extreme example of this is provided by
deterministic interactions ($M=0$ or 1), such as those defining
fixed-routing loss networks and some statistical mechanical models
(see examples in Section \ref{ssec.ex}).  In these cases, the
configurations resulting from forbidden cases ---that is those of the
form $\underline\xi+\delta_\theta$ with $M(\theta\vert \underline
\xi)=0$--- are \emph{not} acceptable as initial configurations for the
interacting process (thus, they will not be generated by it).  The
remaining configurations will be referred to as \emph{acceptable
  configurations}.

\paragraph{Further notation}
The \emph{support} of a function $F$ on $\N^{\ani}$ is the
set
\begin{equation}
  \label{eq:r10}
\supp(F) \;=\; \Bigl\{x\in \Z^d : \,\exists  \,
\ga\ni x \,,\, \underline\eta\in \N^{\ani} \hbox{ such that }
F(\underline\eta) \neq F(\underline\eta + \delta_\gamma)
\Bigr\}\;.
\end{equation}

In this paper we reserve the symbol $\PP$ for the probability measure
on the random environment, while the combination $\Q^{\jota}$ (perhaps
with further embellishments) will denote probability with respect to
processes defined by a fixed environment $\jota$.  We shall use
superindices for the time coordinate of space-time animal
configurations:
$\underline\eta_\Lambda^t\in\N^{\ani_\Lambda}\times\R$.  For finite
sets $A$, the symbol $|A|$ will denote the cardinality of $A$.  For
animals, the notation $|\gamma|$ will mean $|V(\gamma)|$.  We shall
use a capital letter to denote a space-time point. As a default, we
shall use the corresponding lowercase letter to denote its space
component, and a subscript to identify its time component, ex.\
$X=(x,t_X)\in\Z^d\times \R$.

\subsection{Examples}\label{ssec.ex}

Some of the processes fitting our framework are the following.

\subsubsection*{Point processes}

In these models the animals are subsets of $\Z^d$, hence
$V(\ga)=\ga$.  A large variety of point processes has been introduced
in the literature.  Our framework applies to their space-discretized
version.  The \emph{area-interaction point processes} of Baddeley and
van Lieshout (1995\nocite{badlie95}) involve identical animals defined
by the translations of a fixed compact set $G\subset\Z^d$.  An animal
of the form $\gamma = x + G$ is called a \emph{grain} of \emph{germ}
$x$.  Animal configurations are labelled by the corresponding
germs through the identification $\ani\supset\{x+G: x\in A\}
\longleftrightarrow A\subset\Z^d$.  With this identification,
the interaction function takes the form
\begin{equation}
  \label{eq:1}
M(x\vert A) \;=\; F\Bigl[\Bigl|(x+G)\cap (A\oplus G)
\Bigr|\Bigr]
\end{equation}
for some $F:\N\rightarrow [0,1]$, where $A \oplus G = \bigcup_{y \in
  A}\, \{y + G\}$.  The process is attractive if $F$ is increasing and
repulsive otherwise.  The former case corresponds to the
\emph{penetrable sphere model} introduced by Widom and Rowlison
(1970)\nocite{widrow70}.

In addition, the model is specified by internal-clock rates $w(x)$.
In the non-random version they are usually independent of $x$.  For
the disordered version, the independence hypotheses \reff{eq:10} leads
to \emph{site} disorder:
\begin{equation}
  \label{eq:6}
\jota \;=\; \{J_x : x\in \Z^d\}
\end{equation}
for independent random variables $\{J_x\}$.  The internal clock rates
should be of the form
\begin{equation}
  \label{eq:7}
w^\jota(x) \;= \;w(\{J_y: y\in x+G\})\;.
\end{equation}

As further examples of point processes within our framework we
mention: (i) the \emph{Strauss process} (Strauss, 1995\nocite{str75}),
which does not involve grains and the interaction function depends on
the number of pairs of points closer than a fixed threshold $r$, and
(ii) the \emph{perimeter-interaction process} where the interaction
matrix depends on the overlapping of perimeters of grains.  The latter
is a particular instance of the generalization propossed by Baddeley,
Kendall and van Lieshout (1996)\nocite{bkl96}.  Our results apply to
the corresponding lattice versions with site disorder.

\subsubsection*{Fixed-routing loss networks}

These models are defined on a graph $\mathcal{G}$ whose links define
possible calls or connections.  An example is $\Z^d$ with the usual
nearest-neighbor links, or with a more general family of links
$(x,y)\in\Z^d\times \Z^d$ with vertices $x$ and $y$ not necessarily at
distance one.  An animal $\gamma$ is a connected finite subgraph
possibly subjected to some further restriction, depending on the
model (eg.\ not having a loop, or forming a closed circuit).  The
interaction matrix is often deterministic, preventing calls that would
cause a link to be used beyond a predetermined maximal capacity.  More
generally, the interaction matrix embodies some penalization scheme
for the multiple use of links.

The non-random rates $w(\ga)$ usually decrease with the number of
links in $\ga$.  A general type of disorder satisfying the
independence condition \reff{eq:10} is \emph{site-link} disorder:
\begin{equation}
  \label{eq:5}
  \jota \;=\; \{J_x : x\in V(\mathcal{G})\} \,\bigcup\,
\{J_{x,y} : (x,y)\in L(\mathcal{G})\}
\end{equation}
with all the random variables $J_x$, $J_{x,y}$ being independent.  By
$V(\mathcal{G})$ and $L(\mathcal{G})$ we denote, respectively, the set
of vertices and of links of the graph $\mathcal{G}$.  Condition
\reff{eq:10} is satisfied by rates of the form
\begin{equation}
  \label{eq:8}
  w^\jota(\gamma) \;=\; w(\{J_x: x\in V(\ga)\}\,,\, \{J_{x,y} :
  (x,y)\in L(\ga) \})\;.
\end{equation}
We emphasize that, due to our hypothesis \reff{eq:3}, our disordered
loss network must involve connections not exceeding a radius $\ell_1$.

\subsubsection*{Models motivated by statistical mechanics}

These are models of random geometrical objects coming from the study
of spin systems in statistical mechanics. The associated
birth-and-death processes have the statistical mechanical measure as
invariant measure, and can be seen both as a convenient tool to study
such measure and as a feasible simulation algorithm for it. We
mention two of these models.  Due to requirement \reff{eq:3}, the
disordered version considered here correspond to a ``chopped'' variant
of the corresponding models.

\paragraph{\emph {Random cluster model}}

This model plays an essential role in the study of the Potts model.
We refer the reader to Grimmett (1995)\nocite{gri95} for its detailed
study as well as for references to the original articles.  The
animals, called \emph{clusters}, are all finite subgraphs of a given
graph $\mathcal{G}$ (whose vertices are usually the sites of $\Z^d$,
but whose links may include non-nearest-neighbor pairs).  The
interaction matrix forbides the appearance of two clusters having a
common vertix.  Furthermore, the process has internal-clock rates of
the form
\begin{equation}
  \label{wei}
  w(\ga) \;=\;
\prod_{(x,y)\in L(\ga)}  \Bigl({J_{x,y}\over 1-J_{x,y}}\Bigr)\;
\prod_{x\in V(\ga)} \Bigl({1\over J_x}\Bigr)
\end{equation}
for appropiate parameters $J_{x,y}$ and $J_x$ [often denoted
$p(x,y)$ and $q(x)$].  A disordered version, satisfying the
independence condition \reff{eq:10}, is obtained by turning these
parameters into independent random variables [site-link disorder
\reff{eq:5}].

\paragraph{\emph{Ising contour model}}

This is the model resulting from mapping into Peierls contours the
typical spin configurations of one of the low-temperature pure phases
of the Ising model.  Contours are hypersurfaces formed by
$(d-1)$-dimensional unit cubes centered at points of $\Z^d$ and
perpendicular to the edges of the dual latice $\Z^d+(\frac{1}{2},
\cdots, \frac{1}{2})$.  These cubes are called \emph{plaquettes}.  Two
plaquettes are \emph{adjacent} if they share a $(d-2)$-dimensional
face.  A contour is a finite family of plaquettes such that (i) the
family can not be partitioned into two subfamilies with no adjacent
plaquettes, and (ii) every $(d-2)$-dimensional face is covered by an
even number of plaquettes of the family. Geometrically, a contour
corresponds to a connected closed (hyper)surface.  The set of vertices
$V(\ga)$ is the set of centers of the plaquettes forming $\ga$.

The contour-model equilibrium measure is the invariant measure of a
birth-and-death process with a deterministic interaction matrix which
prevents the appearance of two contours with adjacent plaquettes. The
clock rates take the form
\begin{equation}
\label{newpesos}
w(\ga) = e^{-\beta \sum_{x\in V(\ga)} J(x)},
\end{equation}
where the parameter $\beta$ is interpreted as inverse temperature and
each $J(x)$ is a \emph{coupling constant} associated to the plaquette
centered at $x$.  In the disordered version, these constants are
independent random variables [site disorder \reff{eq:7}].

\subsection{Main result}

\begin{teo}
  \label{155}
Consider an animal model as defined above, in particular satisfying
  \reff{eq:3}, \reff{eq:11} and \reff{eq:10}.
  Assume the disorder satisfies
\begin{itemize}
\item[(i)]
  \begin{equation}
    \label{eq:r1}\label{pf2}
\que \;:=\;
\E \biggl[ \ln^{\textstyle a} \Bigl(1 + \sup_x \sum_{\gamma \ni x}
w^\jota(\gamma) \Bigr) \biggr] \;<\; \infty
  \end{equation}
  for some
\begin{equation}
  \label{eq:r2}\label{pf1}
  a \;>\; 2 d^2 \Bigl( 1 + \sqrt{1+ {1 \over d}} + {1\over
  2d}\Bigr)\;.
\end{equation}
\end{itemize}
Then there exists an $\varepsilon=\varepsilon(d,a,\que)>0$, with a
monotonically decreasing dependence $\que\to\varepsilon(d,a,\que)$,
such that if there exists a function $S:\ani\to [1,\infty)$ with
\begin{itemize}
\item[(ii)]
  \begin{equation}
    \label{eq:r3}
    \E\biggl[\, \sup_{\gamma \in \ani}
\frac{1}{S(\gamma)} \sum_{\theta: \theta
  \not\sim \gamma} S(\theta) \,  w^{\jota}(\theta) \biggr]
\;\le\; \varepsilon\;,
  \end{equation}
\end{itemize}
then there exist constants $m>0$, $q>0$ such that for $\P$-almost all
configurations ${\bf J}$ the following is true.
There exist constants
$C_x^\jota\ge 0$, $x\in\Z^d$, such that:

\begin{enumerate}

\item {\bf Existence and uniqueness.}  For each
  $\Lambda\subset\Z^d$ there exists a unique time-invariant
  process $\qint^\jota_{\Lambda}$ on $\N^{\G_\Lambda}\times \R$ with
  generator
  \reff{gen}. The process has a unique invariant
  measure $\mu_\Lambda^\jota$.

\item \label{p2} {\bf Time convergence.} For each $\Lambda\subset\Z^d$
  and each acceptable configuration $\underline \xi_\Lambda
  \in\N^{\ani_\Lambda}$, there exists a unique process
  $\qint^\jota(\,\cdot\,| \,\underline\xi_\Lambda)$ on
  $\N^{\G_\Lambda}\times \R_+$ defined by the generator \reff{gen} and
  the initial configuration $\underline \xi_\Lambda$.  Furthermore,
  the process converges at a superpolynomial rate to the measure
  $\mu_\Lambda^\jota$ on $\N^{\G_\Lambda}$.  Explicitly, the following
  inequality holds for any function
$f$ on $\N^{\ani_\Lambda}$:
    \begin{equation}
      \label{334}
      \sup_{\underline \xi_\Lambda} \left\vert \mu^\jota_\Lambda f -
      \qint^{\jota}( \,f(\underline\eta_\Lambda^t)\,|
      \,\underline\xi_\Lambda) \right\vert \;\le\;\|f\|_\infty \,
C^\jota(\supp(f)) \, \exp\Bigl[-m\,\ln^q(1+t)\Bigr]
    \end{equation}
    with $C^\jota(\supp(f))=\sum_{x\in\supp(f)} C_x^\jota$.  The
    supremum is taken over all acceptable configurations in
    $\N^{\ani_\Lambda}$.

  \item {\bf Space convergence.}  As $\Lambda\to\Z^d$,
    $\mu^\jota_\Lambda$ converges weakly to
    $\mu^\jota:=\mu^\jota_{\Z^d}$ at an exponential rate.  More
    precisely, for any function $f$ depending only on animals
    contained in a finite set $\Lambda$:
\begin{eqnarray}
  \Bigl|\mu^\jota f- \mu^\jota_\Lambda f\Bigr|\;\le \;
\|f\|_\infty\sum_{x\in\supp(f)} C_x^\jota \, \exp\Bigl[ -m\, {\rm
  dist}(x,\Lambda^{\rm c})\Bigr]\;.
\label{73a}
\end{eqnarray}

\item \label{p4}{\bf Exponential mixing.} For any functions $f$ and
  $g$ depending on animals contained in an arbitrary set
  $\Lambda\subset\Z^d$:
    \begin{equation}
      \label{101}
      \Bigl\vert \mu^\jota_\Lambda (f g) - \mu^\jota_\Lambda f\,
\mu^\jota_\Lambda g  \Bigr\vert \;\le\;
      \|f\|_\infty\,\|g\|_\infty\,\sum_{\scriptstyle x\in
        \supp(f),\atop\scriptstyle y\in \supp(g)}
C_x^\jota \, C_y^\jota  \exp\Bigl[-m\,\|x-y\|/2\Bigr]\;.
    \end{equation}

  \item {\bf Perfect simulation.} The invariant measure $\mu^\jota$
    (or any of the $\mu^\jota_\Lambda$) can be perfectly simulated
    using the ancestors algorithm introduced by Fern\'{a}ndez, Ferrari
    and Garcia (2002\nocite{ffg00}).
\end{enumerate}

\end{teo}

These properties
are similar to those of the deterministic process
obtained by a fixed translation-invariant asignment of $\jota$.  The
major difference is the possible slowing-down of the time convergence
due to the subexponential time dependence in the right-hand side of
\reff{334}.

We remark that while the size function $S$ involved in hypothesis
\reff{eq:r3} is completely arbitrary, it is clear that it conveys some
idea of the ``mass'' or ``might'' (Dobrushin's, 1996\nocite{dob96},
terminology) of animals, in the sense that the larger $S(\gamma)$ the
larger
the set of animals
incompatible with $\gamma$.  In fact, the
natural choice for the examples presented in Section \ref{ssec.ex} is
the size of the relevant set of vertices of each animal (area of the
grain, length of the call, perimeter of the contour, etc.). But, of
course, ``natural'' is not a synonym of ``optimal''.

In most cases of interest, condition \reff{eq:r3} is morally more
limitant than \reff{pf2}.  To formalize this fact, let us denote
\begin{equation}
    \label{eq:r30-1}
    \Upsilon^\jota \;=\; \sup_x \sum_{\gamma \ni x}
w^\jota(\gamma)
\end{equation}
and
\begin{equation}
    \label{eq:r30}
    \Psi^\jota \;=\; \sup_{\gamma \in \ani} \,
\frac{1}{S(\gamma)} \sum_{\theta: \theta
  \not\sim \gamma} S(\theta) \,  w^{\jota}(\theta)\;,
\end{equation}
and let us introduce the \emph{halo} of an animal, namely the region
of $\Z^d$ around it beyond which no incompatible is
possible:
\begin{equation}
  \label{eq:rr.50}
  H(\gamma) \;:=\; \bigcap \,\Bigr\{W\subset\Z^d : V(\theta)\subset W
\Longrightarrow \theta\sim\gamma \Bigr\}^c\;.
\end{equation}
This halo is a finite set because of assumption \reff{eq:11}.  In the
examples of Section \ref{ssec.ex} $H(\gamma)=V(\gamma)$.  In addition
we denote
\begin{equation}
  \label{eq:rr.52}
u_1 \;=\; \inf_{\gamma\in\ani:\gamma\neq\emptyset}
\frac{|H(\gamma)|}{S(\gamma)}
\quad,\quad
u_2 \;=\; \sup_{\gamma\in\ani} \frac{|H(\gamma)|}{S(\gamma)}
\end{equation}
and define
\begin{equation}
  \label{eq:rr.53}
\Xi^\jota \;=\;
 \sup_x \sum_{\theta \ni x} \left|H(\theta)\right| \,w^\jota(\theta)\;.
\end{equation}
In terms of these quantities we have the bounds
\begin{equation}
  \label{eq:rr.51}
  \Upsilon^\jota\;\le\; \Xi^\jota \qquad\hbox{and}\qquad
\Psi^\jota \;\le\; \frac{u_2}{u_1}\,\Xi^\jota
\end{equation}
which lead to the following corollary of our main result.
\begin{coro}\label{coro.main}
  Consider an animal model as defined in Section \ref{s.basic} such
  that
  \begin{equation}
    \label{eq:rr.53.1}
u_1>0 \quad\mbox{and}\quad u_2<\infty    \;.
  \end{equation}
  Then, there exists $\widetilde\varepsilon>0$ such that if
  $\E(\Xi^\jota)\le \widetilde\varepsilon$ there exist constants $m,
  q>0$ such that properties 1 to 5 of Theorem \ref{155} hold for
  $\P$-almost all disorder configurations ${\bf J}$.
\end{coro}

The conditions \reff{eq:rr.53.1} are trivially satisfied by all the
examples in Section \ref{ssec.ex}, for which, in fact, $u_1=u_2=1$.
In our bounded-diameter setting, they are also satisfied if the
function size is chosen translation invariant.
\smallskip

\proof We first observe that for any $a\ge 1$ there exists $C_a$ such that
$\ln^a(1+x)\le C_a \,x$ for all $x\ge 0$.  Therefore,
\begin{equation}
  \label{eq:rr.51.1}
  \E\Bigl[\ln^{\textstyle a}(1+\Upsilon^\jota)\Bigr]
\;\le\;  C_a\, \E\Bigl[(\Upsilon^\jota)\Bigr] \;.
\end{equation}
Pick an $a$ as in \reff{pf1} and take
\begin{equation}
  \label{eq:rr.54}
  \widetilde\varepsilon \;=\; \min\Bigl\{\frac{1}{C_a}\,,\,
\frac{u_1}{u_2} \,\varepsilon(d,a,1)\Bigr\}\;,
\end{equation}
where $\varepsilon(d,a,\que)$ is the function given in Theorem
\ref{155}.  From \reff{eq:rr.51.1} and the leftmost inequality in
 \reff{eq:rr.51} we obtain
 \begin{equation}
   \label{eq:rj.1}
   \que \;=\;  \E\Bigl[\ln^{\textstyle a}(1+\Upsilon^\jota)\Bigr]
\;\le\; C_a\, \widetilde\varepsilon \;\le\; 1\;.
 \end{equation}
 From the rightmost inequality in \reff{eq:rr.51} and the
 monotonicity of the function $\que\to\varepsilon(d,a,\que)$ we
 obtain
\begin{equation}
  \label{eq:rr.56}
\E[\Psi^\jota] \;\le\; \frac{u_2}{u_1}\,\widetilde\varepsilon\;
\;\le\; \varepsilon(d,a,1) \;\le\; \varepsilon(d,a,\que)\;.
\end{equation}
Thus, both hypothesis \reff{pf2} and \reff{eq:r3} of
Theorem \ref{155} are satisfied.  $\square$

\section{The basic construction and the key lemma}
\label{const}

The proof of Theorem \reff{155} relies on two ingredients: a graphical
construction of birth-and-death processes propossed by Fern\'{a}ndez,
Ferrari and Garcia (1998\nocite{ffg98}, 2001\nocite{ffg01}), and a multiscale
argument following Klein (1994)\nocite{kle94}. In this section we
review the graphical construction and state the key result obtained
from Klein's argument.

\subsection{Graphical construction}

The construction is performed for each \emph{fixed} environment
$\jota$.

\paragraph{The free process}
We associate to each animal $\gamma$ an independent marked Poisson
process $N_{\gamma}$ with rate $w^{\jota}(\ga)$. We call
$T_{k}(\gamma)$, $\gamma \in \ani$, the ordered time-events of
$N_{\gamma}$ with the convention that $T_0 (\gamma) < 0 < T_1
(\gamma)$. For each occurrence time $T_i (\gamma)$ of the process
$N_{\gamma}$ we choose independent marks $S_i (\gamma)$ exponentially
distributed with mean $1$ and $Z_i(\gamma)$ uniformly distributed in
$[0,1]$.  The \emph{free animal process} is the process in which at
each Poisson time-event $T_i (\gamma)$ (a copy of) the animal $\gamma$
appears and lasts $S_i (\gamma)$ time units. It is convenient to
identify each marked point
$(\gamma, T_k (\gamma), S_k (\gamma),
Z_k(\gamma))$ with $(\gamma \times [T_k(\gamma), T_k(\gamma) +
S_k(\gamma)]\,,\, Z_k(\gamma))$, the {\em cylinder} with {\em basis}
$\gamma$, {\em
  birth time} $T_k(\gamma)$ {\em lifetime} $S_k(\gamma)$ and
\emph{mark} $Z_k(\gamma)$. This identification turns the free process
into a measure on the space of cylinders. Denoting
$C=(\gamma,t,s,z)$, we use the notation
\[
\basis(C) = \gamma, \,\,\, \birth(C)=t, \,\,\, \life(C)= [t,t+s],
\,\,\, \mark(C)= z
\]
The marks will be used later on to define the (interacting)
birth-and-death process with generator \reff{gen}.  Let us denote
$\Q^\jota$ the probability measure on $\Z^d\times \R$ corresponding to
the free process (it always exist, being a countable product of marked
Poisson processes).

\paragraph{Backwards oriented percolation}
We first extend our definition of incompatibility to cylinders in the
natural way: Two cylinders $C$ and $C'$ are incompatible if they have
incompatible bases and they are simultaneously alive at some instant
of time.  That is, $C \not\sim C'$ if and only if $\basis(C) \not\sim
\basis(C')$ and $\life(C) \cap \life(C') \neq \emptyset$; otherwise $C
\sim C'$.

Let us now fix a family of cylinders $\rea$ (for instance, obtained as
a realization of the free process of the previous paragraph).  We
define the \emph{ancestors} of a cylinder $C\in\rea$ as the set
\begin{equation}
{\bf A}_1^{C} = \Bigl\{C' \in \rea: C' \not \sim C \mbox{ and }
\birth(C') \leq \birth(C) \Bigr\}\;.
\end{equation}
Recursively for $n \geq 1$, the $n$th generation of ancestors of $C$
(in $\rea$) is
\begin{equation}
{\bf A}_n^{C} = \Bigl\{ C'': C'' \in {\bf A}_1^{C'} \mbox{ for some } C'
  \in {\bf A}_{n-1}^{C}\Bigr\}\;.
\end{equation}
[In fact, ${\bf A}_n^{C}={\bf A}_n^{C}(\rea)$, but we shall omit this
dependence except if it is crucially needed.]  Likewise, for an
arbitrary space-time point $(x,t) \in {\Z}^d \times {\R}$ we define
its set of ancestors as the set of cylinders that contain it
\begin{equation}
{\bf A}_1^{(x,t)} = \Bigl\{C \in \rea: \basis(C) \ni x, \life(C) \ni
t\Bigr\}\;,
\end{equation}
and, recursively,
\begin{equation}
\label{rp3}
{\bf A}_n^{(x,t)} = \Bigl\{ C'': C'' \in {\bf A}_1^{C'} \mbox{ for some }
C' \in {\bf A}_{n-1}^{(x,t)}\Bigr\}\;.
\end{equation}

\begin{defin} We say that \emph{there exists backwards oriented
    percolation} in $\rea$ if there exists a space-time point $(x,t)$
  such that ${\bf A}_n^{(x,t)} \neq \emptyset$ for all $n$, that is,
  there exists a point with infinitely many generations of ancestors.
\end{defin}

The {\em clan} of a space-time point $(x,t)$, resp.\ of a cylinder
$C$, is the union of the corresponding ancestors:
\begin{equation}\label{eq:nue.1}
{\bf A}^{(x,t)} = \bigcup_{n \geq 1} {\bf A}_n^{(x,t)} \qquad ,
\qquad {\bf A}^{C} = \bigcup_{n \geq 1} {\bf A}_n^{C}\;.
\end{equation}

Other quantities that will be used later are the time-length,
space-diameter and space-size of the clan of a point $(x,t)$:
\begin{eqnarray}
\label{teele}
\tl({\bf A}^{(x,t)}) &=& t - \inf\Bigl\{ s: \life(C) \ni s, \mbox{ for
  some } C \in {\bf A}^{(x,t)}\Bigr\}\\
\label{esedobleve}
\sd({\bf A}^{(x,t)}) &=& {\rm diam}\Bigl(\bigcup_{C \in {\bf A}^{(x,t)}}
\basis(C)\Bigr)\\
\ss({\bf A}^{(x,t)}) &=& \Bigl|\bigcup_{C \in {\bf A}^{(x,t)}} \;
\basis(C)\Bigr|\;.
\end{eqnarray}

\paragraph{The interacting birth-and-death process}
If
\begin{equation}
  \label{eq:4}
  \Q^\jota(\{\hbox{no backwards percolation}\}) \;=\; 1\;,
\end{equation}
the stationary process with generator \reff{gen} can be constructed by
``cleaning'' the free process defined above.  For completeness we
present a summary of this construction here. The reader is referred to
Fern\'{a}ndez, Ferrari and Garcia (2001) for details.  The idea is to
start from first ancestors (``Eves'') and classify cylinders into
\emph{kept} or \emph{erased} according to the test determined by the
interaction function.  Cylinders that are born in presence of a kept
ancestor and that fail the test are erased, all the others are kept.

Explicitly, let $\rea$ be a cylinder configuration in the set
\begin{equation}
  \label{eq:9}
\C=\{\rea \hbox{ without backwards percolation} \}  \;.
\end{equation}
Since all clans in $\rea$ are finite, each cylinder has a well
defined, finite number of ancestors. Therefore the configuration can
be decompossed in the form $\rea=\cup_{n\ge 0} \rea_n$, where
$\rea_n:= \{C\in\rea: \A_{n}^C \neq \emptyset,\;\A_{n+1}^C
=\emptyset\} $.  The sets
$\kep$ and $\era$ of kept and erased cylinders
are defined inductively as follows.  Starting with
$\kep_0=\rea_0$ and
$\era_0=\emptyset$ (cylinders without ancestors are kept, as they do not
need to pass any test), we define, recursively,
\begin{eqnarray}
  \label{202}
  \kep_n &=& \Bigl\{C\in \rea_n\setminus
\cup_{i=0}^{n-1}(\era_i\cup\kep_i)\,:\, \mark(C)\le \widetilde M(C|
\cup_{i=0}^{n-1}\kep_i)\Bigr\} \nonumber\\
\ \\
\era_n &=&\rea_n\setminus[\kep_n\cup\,
\cup_{i=0}^{n-1}(\era_i\cup\kep_i)]\nonumber
\end{eqnarray}
where $\widetilde M(C|\kep') = M(\basis (C)|\{\basis(C')\,:\,C'\in \kep',
\birth(C)\in \life (C')\})$ .
We denote the set of kept cylinders as
$\kep(\rea)=\cup_n \kep_n$ and the
set of erased cylinders as $\era(\rea)=\cup_n\era_n$. Clearly
\begin{itemize}
\item[(C1)] $\rea$ is the disjoint union of $\kep(\rea)$ and
  $\era(\rea)$, and
\item[(C2)] The event $\{C\in\kep(\rea)\}$ is measurable with respect
  to the sigma field generated by $\A^C$, in fact $\kep(\A^C) =
  \kep(\rea)\cap\A^C$.  In words, it is sufficient to know the
  (finite) clan of $C$ to know whether $C$ is kept or erased.
\end{itemize}

The stationary animal process $\overline\Q^\jota$ is defined by the
sections of the kept cylinders:
\begin{equation}
  \label{152} \eta^{\jota,t}(\ga,\rea) \; = \;\sum_{C\in\kep(\rea)} \one
  \Bigl\{\basis(C)=\ga,\, \life(C)\ni t\Bigr\}\;.
\end{equation}
If $\Q^\jota(\C)=1$, the process $\overline\Q^\jota$ is Markovian and
has generator \reff{gen}; that is,
\begin{equation}
  \label{206}
  {d \over dt} \overline\E^\jota f(\eta^t(\rea)) \;=\;
\overline\E^\jota\, A_\Lambda^\jota f(\eta^t(\rea))
\end{equation}
This fact is proven in Fern\'{a}ndez, Ferrari and Garcia (2001) for the
homogeneous case.  This proof extends to the inhomogeneous case in an
obvious manner. Let us denote $\mu^\jota$ the distribution of any
$t$-section $\underline\eta^{\jota,t}$, which, by construction, is
independent of $t$.
We shall determine the properties of $\mu^\jota$ by studying the law of
$\eta^{\jota,0}$, the stationary (interacting) birth-and-death process
at time zero.

As observed in (C2) above, the presence or absence of contours
intersecting a region $V\subset \Z^d$ at time $t$ depends only
on the clans of the cylinders alive at time $t$ whose bases intersect
$V$, that is, on
\begin{equation}
  \label{eq:12}
\A^{V,t}\;:=\; \Bigl\{C'\in\A^C\,:\, \basis(C)\cap
V \neq\emptyset,\, \life(C)\ni t\Bigr\} \;.
\end{equation}
In particular the function $\eta^{\jota,t}(\ga,\,\cdot\,)$ defined by
\reff{152} is in fact a (deterministic) function only of $\A^{\ga,t}$
[$=\A^{V(\ga),t}$].  More precisely, if $\rea$ and $\rea'$ are two
cylinder configurations such that
$\A^{\gamma,t}(\rea)=\A^{\gamma,t}(\rea')$, then
$\eta^{\jota,t}(\ga,\rea)=\eta^{\jota,t}(\ga,\rea')$.  We code this
fact as the identity (slightly abusive from the notational point of
view):
\begin{equation}
  \label{209}
 \eta^{\jota,t}(\ga,\rea)=\eta^{\jota,t}(\ga,\A^{\ga,t})\,.
\end{equation}

\subsection{The key lemma}

Let us call a sequence of cylinders $C_1, C_2, \ldots, C_n$ an {\em
  open path} if $C_2 \not\sim C_1 , \ldots ,C_n \not\sim C_{n-1} $ and
$\birth(C_{i+1})<\birth (C_i)$ for all $i$. Given a cylinder
configuration $\rea$ and two space-time points $X=(x,t)$ and $Y=(y,s)$,
with $s \leq t$, we say that $X$ and $Y$ are \emph{connected}
(in the configuration $\rea$) if there exists an open path $C_1,
C_2 ,\ldots, C_n$ such that $x \in
\basis(C_1)$, $t \in \life(C_1)$, $y \in \basis(C_n)$, $s \in
\life(C_n)$.  [Equivalently, $Y$ is in (the interior
of a cylinder belonging to) the clan of ancestors of $X$.]
The existence of such a connection defines an event denoted $X
\rightarrow Y$.

For a given realization of the environment, the {\em connectivity
  function} is defined by
\begin{equation}
G^{\bf J}(X,Y) = {\Q}^{\bf J} \{ X \rightarrow Y \}\;.
\end{equation}

Many of the properties stated in the main theorem are a direct consequence of
the following result

\begin{lema}
\label{main}
Under the hypotheses of Theorem \ref{155} [namely, \reff{eq:3},
\reff{eq:11}, \reff{eq:10} and \reff{eq:r2}], but replacing
\reff{eq:r1} by an inequality
  \begin{equation}
    \label{eq:rj1}
\E \biggl[ \ln^{\textstyle a} \Bigl(1 + \sup_x \sum_{\gamma \ni x}
w^\jota(\gamma) \Bigr) \biggr] \;\le\; \que\;,
  \end{equation}
  there exist $m>0$ and $q_0(a,d)>1$ such that for all $q \in (1,q_0)$
  there is a value $\varepsilon(d,a,m,q,\aleph)>0$ so that the
  validity of \reff{eq:r3} implies that for every $x \in {\Z}^d$
\begin{equation}
\label{conect}
G^{\bf J}\Bigl((x,t),(y,s)\Bigr) \;\leq\; C^{\bf J}_x\, \exp\Bigl\{ -
  m  \max \Bigl[ \norm{ x-y},\, \ln^q (1
    + |t-s|) \Bigr]   \Bigr\}
\end{equation}
for all $y \in {\Z}^d$ and $t,s \in {\R}$, where the constants
$C_x^{\bf J}=C_x^{\bf J}(\ell_1+\ell_2)$ are finite for $\PP$-almost
every environment ${\bf J}$.
\end{lema}

The fact that the inequality \reff{eq:rj1} is used implies that the
function $\que\to\varepsilon(\,\cdots\,,\que)$ can be chosen
to be decreasing, as stated in Theorem \ref{155}.

We will prove this lemma by performing a multiscale analysis similar
to the one used by Klein (1994)\nocite{kle94} in his work on extintion
of contact process in a random environment. That will be done in
Section \ref{s.coar}. We first discuss how Theorem \ref{155} follows
from the bound \reff{conect} for the connectivity function.

\section{How the theorem follows from the key lemma}
\subsection{Relation with percolation properties}

The following theorem
relates properties of the measure $\mu^{\bf J}$
with properties of the percolation model.
\begin{teo}
  \label{225}
  Assume that for a given ${\bf J}$ there is no backwards oriented
  percolation with $\Q^{\bf J}$-probability one.  Then,
\begin{enumerate}
\item {\bf Uniqueness.} The measure $\mu^{\bf J}$ is the unique
  invariant measure for the process $\eta^{\jota,t}$.

\item {\bf Time convergence.} For any function $f$ with finite
  support,
 \begin{eqnarray}
   \lefteqn{ \sup_{\underline\xi_\Lambda} \left\vert
      \mu^\jota_\Lambda f - \overline
      \E^{\jota}\Bigl( \,f(\underline\eta_\Lambda^t)\,\Bigm|
      \,\underline\xi_\Lambda\Bigr) \right\vert} \nonumber\\
   &\le& 2\,\|f\|_\infty\,\sum_{x\in\supp(f)}\Bigl[\Q^\jota\Bigl\{
   \tl(\A^{(x,0)})>bt \Bigr\} +
   e^{-(1-b)t}\, \E^\jota\Bigl(\ss(\A^{(x,0)})\Bigr)\Bigr] \label{335}
   \end{eqnarray}
for any $b\in(0,1)$.

\item {\bf Space convergence.} As $\Lambda\to\Z^d$, $\mu^{\bf
    J}_\Lambda$ converges weakly to $\mu^{\bf J}$. More precisely, if
  $\supp(f)\subset\Lambda$, then
\begin{equation}
  { \Bigl|\mu^\jota f- \mu^\jota_\Lambda f\Bigr|}\;\le
  \;2\,\|f\|_\infty\,\sum_{x\in\supp(f)}\Q^\jota\Bigl\{
\sd(\A^{(x,0)})\ge \dist(\{x\}, \Lambda^c) \Bigr\}\, .\label{73ab}
\end{equation}

\item {\bf Mixing.} For $f$ and $g$ with finite support,
\begin{eqnarray}
\lefteqn{\Bigl\vert \mu^\jota_\Lambda (f g) - \mu^\jota_\Lambda f\,
\mu^\jota_\Lambda g  \Bigr\vert\;\le\;
4\,\|f\|_\infty\,\|g\|_\infty\,\label{p73b}}\nonumber\\
 &&\quad \times\sum_{\scriptstyle
    x\in \supp(f),\atop\scriptstyle y\in \supp(g)}
  \Bigl[\Q^\jota\Bigl\{\sd(\A^{(x,0)})\ge
  \norm{x-y}/2\Bigr\}+\Q^\jota\Bigl\{\sd(\A^{(y,0)})\ge
\norm{x-y}/2\Bigr\}\Bigr]\;.
\nonumber\\
\
\end{eqnarray}
\end{enumerate}
\end{teo}

\proof Items 1, 2 and 3 follow from displays (4.6), (4.7) and (4.9) of
Theorem 4.1 of Fern\'{a}ndez, Ferrari and Garcia (2001). The analogous
of item 4 is stated in that theorem with the extra assumption that
there exists a time $h$ such that there is no space-time percolation
in $(0,h)$. We provide here a proof without this hypothesis.

We consider functions $f$ and $g$ such that $\supp(f)=\{x\}$ and
$\supp(g)=\{y\}$, the general case follows by telescoping. Let's fix a
partition $\{\Gamma, \Gamma'\}$ of $\G_\Lambda$. Below we choose
$\Gamma$ formed by animals ``closer to $x$''.  Let $\underline{ A}$
and $\underline{ B}$ be $\Q^\jota$-independent realizations of
$\{C\in\C\,:\, \basis(C)\in \Gamma\}$ and $\underline{ A}'$ and
$\underline{ B}'$ independent realizations of $\{C\in\C\,:\,
\basis(C)\in \Gamma'\}$.  Then $\underline{ A}\cup \underline{ A}'$,
$\underline{ A}\cup \underline{ B}'$, $\underline{ B}\cup \underline{
A}'$, $\underline{ B}\cup \underline{ B}'$ have the same law as
$\C_\Lambda$ and $\underline{ A}\cup \underline{ B}'$ is independent
of $\underline{ B}\cup \underline{ A}'$.  Let $\A^{(x,t)}(\underline
C)$ be the random variable defined in
\reff{eq:nue.1} and
\[
X(\underline{ A},\underline{ A}') \;=\;
f\Bigl(\eta^{\jota,0}(\,\cdot\,, \underline{A}\cup \underline{
A}')\Bigr) \;=\; f\Bigl(\eta^{\jota,0}(\,\cdot\,,
\A^{(x,0)}(\underline{ A}\cup \underline{ A}'))\Bigr) \;.
\]
Analogously
\[
Y(\underline{ A},\underline{ A}') \;=\;
g\Bigl(\eta^{\jota,0}(\,\cdot\,, \underline{ A}\cup \underline{
A}')\Bigr) \;=\; g\Bigl(\eta^{\jota,0}(\,\cdot\,,
\A^{(x,0)}(\underline{ A}\cup \underline{ A}'))\Bigr) \;.
\]
With these definitions, we obtain
\begin{equation}\label{cam:3}
\mu^\jota_\Lambda (f g) - \mu^\jota_\Lambda f\, \mu^\jota_\Lambda g
\;=\; \widetilde\E \Bigl[X(\underline{ A},\underline{ A}')\,Y(\underline{
A},\underline{ A}') -  X(\underline{ A},\underline{
B}')\,Y(\underline{ B},\underline{ A}')\Bigr]\;,
\end{equation}
where $\widetilde\E$ corresponds to a four-fold product of the measure
$\Q^\jota$.   This expression leads us to the bound
\begin{equation}
  \label{rp2}
\Bigl| \mu^\jota_\Lambda (f g) - \mu^\jota_\Lambda f\,
\mu^\jota_\Lambda g \Bigr|
\;\le \;2\,\|f\|_\infty\,\|g\|_\infty\,\widetilde\P(\mathcal{A}^c)
\end{equation}
with
\begin{equation}
\label{cam.1}
\mathcal{A}= \Bigl\{X(\underline{ A},\underline{ A}')=X(\underline{
A},\underline{ B}') \hbox{ and } Y(\underline{ A},\underline{
A}')=Y(\underline{ B},\underline{ A}') \Bigr\}\;.
\end{equation}

We now choose $\Gamma$ as the set of animals intersecting $\{z\in
\Z^d\,:\, \norm{z-x}\le \norm{z-y}\}$ and $\Gamma'$ as its complement.
Then the event $\mathcal{A}$ is verified whenever the bases of the
cylinders in both $\A^{(x,0)}(\underline{ A}\cup \underline{ A}')$ and
$\A^{(x,0)}(\underline{ A}\cup \underline{ B}')$ are contained in
$\Gamma$ and those of $\A^{(y,0)}(\underline{ A}\cup \underline{ A}')$
and $\A^{(y,0)}(\underline{ B}\cup \underline{ A}' )$ are contained in
$\Gamma'$.  The complement of the intersection of these four events
yields
\[
\widetilde\P(\mathcal A^c) \;\le\; 2\,\Q^\jota\Bigl\{\sd(\A^{(x,0)})\ge
\norm{x-y}/2\Bigr\}\,+\,2\,\Q^\jota\Bigl\{\sd(\A^{(y,0)})\ge
\norm{x-y}/2\Bigr\}\,.
\;\square\]
\medskip

For comparison purposes, let us present an alternative mixing bound.
\begin{prop}
Assume that for a given $\jota$ there is no backwards oriented
percolation with $\Q^\jota$-probability one and consider functions $f$
and $g$ with finite support.  Then,
\begin{eqnarray}\label{cam:2}
\lefteqn{\Bigl\vert \mu^\jota_\Lambda (f g) - \mu^\jota_\Lambda f\,
\mu^\jota_\Lambda g \Bigr\vert\;\le\;
2\,\|f\|_\infty\,\|g\|_\infty\,}\nonumber\\
&& \times \sum_{\scriptstyle x\in
\supp(f),\atop\scriptstyle y\in \supp(g)}
\widehat\Q^\jota\Bigl\{\sd(\A^{(x,0)})+\sd(\A^{(y,0)})\ge
\norm{x-y}\Bigr\}\;,
\end{eqnarray}
where $\widehat\Q^\jota$ is the free process obtained doubling the
birth rates of $\Q^\jota$.
\end{prop}

The bound \reff{cam:2} corresponds to standard high-temperature
results in statistical mechanics (see, for instance, the main result
in Bricmont and Kupiainen, 1996\nocite{brikup96}).  Its proof relies
on the very popular technique of ``duplication of variables''.  In
contrast, \reff{p73b} is proven by ``tetra-plication of variables''.
In our general setting, however, we are able to exploit better
our first bound \reff{p73b}.
\medskip

\proof
As above, it is enough to assume $\supp(f)=\{x\}$ and
$\supp(g)=\{y\}$.  Let $\underline{ C}$ and $\underline{ C}'$ be
$\Q^\jota$-independent realizations of $\C_\Lambda$.  We denote
\[
X(\underline{C}) \;=\;
f\Bigl(\eta^{\jota,0}(\,\cdot\,, \underline{C})\Bigr)
\quad\hbox{ and } \quad
Y(\underline{C}) \;=\;
g\Bigl(\eta^{\jota,0}(\,\cdot\,, \underline{C})\Bigr)\;.
\]
The duplication-of-variables identity is:
\begin{equation}\label{cam:4}
\mu^\jota_\Lambda (f g) - \mu^\jota_\Lambda f\, \mu^\jota_\Lambda g
\;=\; \frac{1}{2}\,
\widetilde\E \Bigl[\Bigl(X(\underline{C}) -X(\underline{C}')\Bigr)
\Bigl(Y(\underline{C})-Y(\underline{C}')\Bigr)\Bigr]\;,
\end{equation}
where $\widetilde\E$ corresponds to the measure
$\Q^\jota\times\Q^\jota$.  Let us now consider the event
\begin{equation}\label{cam:6}
\mathcal{B} \;=\; \Bigl\{
\basis\Bigl(\A^{(x,0)}(\underline C\cup \underline C')\Bigr)
\,\bigcap\,
\basis\Bigl(\A^{(y,0)}(\underline C\cup \underline C')\Bigr)
 \,=\, \emptyset \Bigr\}\;,
\end{equation}
and the transformation $T:(\underline C,\underline C') \to
(T\underline C,T\underline C')$, that interchanges the
$(x,0)$-ancestors in $\underline C$ and $\underline C'$:
\begin{eqnarray}
T \underline C &=& \Bigl(\underline C\setminus \A^{(x,0)}(\underline
C)\Bigr) \cup \A^{(x,0)}(\underline C')\nonumber\\
T \underline C' &=& \Bigl(\underline C'\setminus \A^{(x,0)}(\underline
C')\Bigr) \cup \A^{(x,0)}(\underline C)\;.
\end{eqnarray}
Conditioned to $\mathcal{B}$ being true, the distribution of
$(\underline C,\underline C')$ coincides with that of $(T\underline
C,T\underline C')$ (the processes inside and outside each realization
of $\basis\Bigl(\A^{(x,0)}(\underline C\cup \underline C')\Bigr)$ are
independent).  Furthermore, the event $\mathcal{B}$ is $T$-invariant.
Hence,
\begin{equation}\label{cam:7}
\widetilde\E\Bigl[ F(\underline C,\underline C')
\,  \one\{\mathcal{B}\}\Bigr] \;=\;
\widetilde\E\Bigl[ F(T\underline C,T\underline C')
\;  \one\{\mathcal{B}\}\Bigr]
\end{equation}
for each local function $F$ on $\G_\Lambda\times \G_\Lambda$.
But, in the presence of $\mathcal{B}$, the function $F$ involved
in \reff{cam:4} is odd under this transformation.  We conclude that
\begin{equation}\label{cam:10}
\widetilde\E \Bigl[\Bigl(X(\underline{C}) -X(\underline{C}')\Bigr)
\Bigl(Y(\underline{C})-Y(\underline{C}')\Bigr)\,
\one\{\mathcal{B}\} \Bigr]\;=\; 0\;,
\end{equation}
which, by \reff{cam:4}, implies
\begin{equation}
  \label{cam:9}
\Bigl| \mu^\jota_\Lambda (f g) - \mu^\jota_\Lambda f\,
\mu^\jota_\Lambda g \Bigr|
\;\le \;2\,\|f\|_\infty\,\|g\|_\infty\,\widetilde\P(\mathcal{B}^c)\;.
\end{equation}
The proof follows from the observation that
\[
\widetilde\P(\mathcal B^c) \;\le\;
\widehat\Q^\jota\Bigl\{\sd(\A^{(x,0)})+\sd(\A^{(y,0)})\ge
\norm{x-y}\Bigr\} \;.\quad  \square
\]

\subsection{Bounds for the size of the clan}
The
 key inequality \reff{conect} leads to the following bounds for the
probabilities of the time-length and space diameter and size of the
clan of a space-time point.

\begin{prop} Under the hypotheses of Theorem \ref{155}
  there exist $\overline{m}, \widetilde m>0$ and constants $\overline
  C_x^\jota$, $\widetilde C_x^\jota$, and $\widehat C_x^\jota$,
  $x\in\Z^d$, which are finite for $\PP$-almost every environment,
  such that
\begin{eqnarray}
\label{gardel1}
\Q^\jota\Bigl\{\tl({\bf A}^{(x,0)}) > T \Bigr\} &\leq&
\overline C_x^\jota\, \exp\Bigl\{- \overline{m} \,\ln^q (1 + T) \Bigr\}\\[7pt]
\label{gardel3}
\Q^{\jota}\Bigl\{\sd({\bf A}^{(x,0)}) > L \Bigr\} &\le&
\widetilde C_x^\jota\,\exp\{- \widetilde m L \}\\[7pt]
\label{gardel2}
\E^{\jota}\Bigl\{\ss({\bf A}^{(x,0)})\Bigr\} &\leq&
\widehat C_x^\jota\;.
\end{eqnarray}
\end{prop}

\proof
As $\max\,(a,b)\ge (a+b)/2$, inequality \reff{conect} leads to
\begin{equation}
\label{conectbis}
G^{\bf J}\Bigl((x,t),(y,s)\Bigr) \;\leq\; C^{\bf J}_x\, \exp\Bigl\{ -
\frac{m}{2}\, \Bigl[ \norm{x-y} + \ln^q (1
    + |t-s|) \Bigr]   \Bigr\} \;.
\end{equation}
Inequality \reff{gardel1} is then a straightforward consequence of the
fact that
\begin{equation}
\label{eq:13}
  \Q^{\jota}\Bigl\{\tl({\bf A}^{(x,0)}) > T \Bigr\}
\;\le\; \sum_y G^{\jota}\Bigl((x,0) \rightarrow (y,-T)\Bigr) \;,
\end{equation}
while \reff{gardel3} follows from
\begin{equation}
  \label{eq:14}
  \Q^{\jota}\Bigl\{\sd({\bf A}^{(x,0)}) > L \Bigr\}  \;\leq\;
\sum_{y: \norm{x-y} \ge L}\, \int_0^\infty
G^{\jota}\Bigl((x,0) \rightarrow (y,-t)\Bigr) dt \;.
\end{equation}
To obtain \reff{gardel2} we use \reff{gardel3} to bound the right-hand
side of the inequality
\begin{equation}
  \label{eq:rr.45}
\Q^{\jota}\Bigl\{\ss({\bf A}^{(x,0)}) > L \Bigr\}
\;\le\; \Q^{\jota}\Bigl\{\sd({\bf A}^{(x,0)}) > L^{1/d} \Bigr\}
\end{equation}
and sum over $L$.  $\Box$

\subsection{Proof of Theorem \ref{155}}
The bound \reff{conect} provided by the key Lemma \ref{main} implies,
by Borel-Cantelli, the absence of backwards oriented percolation.
Thus, we can apply Theorem \ref{225} for almost-all realizations of
the disorder.
\begin{enumerate}
\item {\bf Existence and uniqueness.\/} Under the hypothesis of no
  percolation in any region $\Lambda$, the process \reff{152} defines
  a stationary process $\eta^t$ in $\Lambda$ by considering those
  cylinders in $\C$ with basis contained in $\Lambda$.  Uniqueness
  follows from the next item.

\item {\bf Time convergence.\/} We use \reff{gardel1} and \reff{gardel2} to
  bound the right hand side of \reff{335} by
\begin{eqnarray}
\label{334a}
\lefteqn{ 2\,\| f \|_\infty\,\sum_{x\in\supp(f)}
\Bigl(\overline C_x^\jota
\exp\Bigl\{- \overline{m} \,\ln^q(1 + bt) \Bigr\}
+ \widehat C_x^\jota \,e^{-(1-b)t} \Bigl)}
   \nonumber \\[7pt]
&&\quad \leq\  2\,\|f\|_\infty\,\sum_{x\in\supp(f)} C_x^\jota  \,
\exp\{-m \ln^q (1 + bt) \},
\end{eqnarray}
for some constant $C_x^\jota < \infty$.

\item {\bf Space convergence.\/} Inequality \reff{73a} follows
  immediately from \reff{73ab} and \reff{gardel3}.

\item {\bf Exponential mixing.\/} The bound \reff{101} is just the
  combination of
\reff{p73b} and \reff{gardel3}.

\item{\bf Perfect simulation.\/} The construction proposed in
  Fern\'{a}ndez, Ferrari and Garcia (2001) works, for almost all disorder
  realizations, due to the absence of backwards oriented percolation.

\end{enumerate}

\section{Proof of Lemma \protect\ref{main}}\label{s.coar}


The lemma is proven by adapting Klein's (1994\nocite{kle94})
multiscale analysis.  The analysis involves several steps that will be
studied separately.

\subsection{General scheme and notation}

The scheme is based on a sequence of linear sizes, called
\emph{scales}, defining boxes of increasing size. Sites are
classified on \emph{regular} or {\em singular} according to the decay
of connectivity functions on a surrounding box.  The definitions must
be tunned up so that the probability for a site to be regular increase
sufficiently fast with the scale.  The presentation of this scheme is
organized as follows:

\begin{enumerate}

\item[(i)] In Section \ref{s.tool} we recopilate a number of
  inequalities needed to prove next inductive step.  In particular,
  the Hammersley-Lieb-Simon inequality is the key to pass from local
  to global decay of the connectivity function. The uniform bound on
  the size of the animals is needed to ensure its validity.

\item[(ii)] In Section \ref{s.multi} we define regularity and
  determine ``good enough'' probabilities for sites to be regular at
  each scale.  Regularity at all scales with such probabilities
  guarantee, by Borel Cantelli, almost sure regularity.

\item[(iii)] Section \ref{s.change} shows the crucial inductive step:
  A ``good enough'' probability of being regular at a given scale
  implies a ``good enough'' probability of being regular at the next
  scale.  This is the heart of the argument.

\item[(iv)] The last step of the proof is the determination of
  conditions so the origin has a ``good enough'' probability of being
  regular at some initial scale.  This is done in Section \ref{s.fin}.

\end{enumerate}

We introduce some notation.  We fix $\ell_0 = \ell_1 + \ell_2 $, the
maximal size of the animals plus the radius of incompatibility -- to
avoid trivialities we assume $\ell_0 > d+1$. For $L>0$ and $x \in
{\Z}^d$, we denote
\begin{equation}
\Lambda[x;L]\;=\; \Bigl\{y \in {\Z}^d: \norm{x-y} \leq L
\Bigr\}
\end{equation}
and define the $\ell_0$-boundary of $\Lambda[x;L]$ as the set
\begin{equation}
\partial_{\ell_0}\Lambda[x;L] \;=\; \Lambda[x; L+ \delta] \setminus
\Lambda[x;L],
\label{r.delta}
\end{equation}
where $\delta>1$ is chosen so the following is true: For any sequence
of animals $\gamma_1, \gamma_2, \ldots , \gamma_n$ with $\gamma_i
\not\sim \gamma_{i+1}$, $i=1,2,\ldots ,n-1$ connecting a point inside
$\Lambda[x;L]$ with a point outside $\Lambda[x; L+ \delta]$, at least
two animals in the sequence are contained in
$\partial_{\ell_0}\Lambda[x;L]$. This condition is necessary to satisfy
the Hammersley-Lieb-Simon (see next).  A possible choice is $\delta =
\delta (\ell_0) = \frac{3(\ell_0 -2)}{2(d-1)} $.

For $X=(x,t) \in {\Z}^d \times {\R}$, $L>0$ and $T>0$ we set
\begin{equation}
\label{boxes}
B_{L,T}(X) \;=\; \Lambda[x;L] \times [t-T, t] \;.
\end{equation}
The vertical, horizontal and complete boundaries of the box
$B_{L,T}(X)$ are defined respectively as:
\begin{eqnarray}
\partial_{V}B_{L,T}(X) &=& \partial_{\ell_0} \Lambda[x;L] \times
[t-T,t],\\
\partial_{H}B_{L,T}(X) &=& \Lambda[x; L+ \delta] \times
\{t-T \},\\
\partial B_{L,T}(X) &=& \partial_{V}B_{L,T}(X) \cup
\partial_{H}B_{L,T}(X).
\end{eqnarray}
Note that, as we are considering backwards oriented percolation, the
face $\Lambda[x; L+ \delta] \times \{t\}$ is excluded from the
boundary.

For any integrable function $H : {\Z}^d \times {\R} \rightarrow {\R}$ we
denote
\begin{equation}
\label{Ulises}
\sum_{Z \in \partial B_{L,T}(X)} \, H(Z) = \sum_{z \in \Lambda[x; L +
  \delta]} \,  H(z, t-T)  + \sum_{z \in \partial_{\ell_0} \Lambda[x;L]} \,
\int_{t-T}^{t}  H(z,s) ds.
\end{equation}

We introduce the notion of connection \emph{within a box}.  Given a
cylinder configuration $\rea$, its \emph{restriction} to a box $B =
\Lambda \times I$, with $\Lambda \subset {\Z}^d$ and $I$ a real
interval, is the family $\underline C_B$ of cylinders obtained by
``restricting'' to $I$ those cylinders of $\underline C$ with bases
inside $\Lambda$.  That is, each $C=(\gamma,t,s,z)\in\underline C$
with $\gamma\in\ani_\Lambda$ defines a cylinder
$C_B=(\gamma,t_I,s_I,z)\in\underline C_B$ with $t_I=\max(t,\inf I)$
and
$s_I= \min(t+s,\sup I)-t_I$.  Given two space-time points, $X=(x,t)$
and $Y=(y,s)$, $X,Y\in B$ with $s \leq t$, the event $X \rightarrow_B
Y$ is formed by all cylinder configurations $\rea$ such that the
configuration $\rea_B$ exhibits an open path connecting $X$ with $Y$.
If this event is true, we say that $X$ and $Y$ are connected \emph{in}
$B$.  The event defines the connectivity function in the region $B$:
\begin{equation}
G_{B}^{\bf J}(X,Y) = {\Q}^{\bf J} \{ X \rightarrow_B Y \}\;.
\end{equation}
As usual, we shall omit the subscript if $B=\Z^d\times\R$.
We write
\begin{equation}
  \label{eq:15}
  G_{B_{L+ \delta, T}(X)}^{\jota}(X, \partial) \;=\;
 \sum_{Z \in \partial B_{L,T}(X)} \,
G_{B_{L+ \delta, T}(X)}^{\jota}(X, Z)\;.
\end{equation}

\subsection{A toolbox:  Inequalities}\label{s.tool}

The first two inequalities we need have been basically proven in the
literature.  The minor adaptations needed for our setting do not
justify a detailed exposition of their proofs.  We content ourselves
with providing appropriate references and indications.

We consider the natural partial order in the space of
cylinder configurations: $ \underline C \leq \underline C'$ if $
\underline C'$ contains all the cylinders in $ \underline C$.  Events
are said \emph{increasing}, resp.\ \emph{decreasing}, if their
characteristic functions are nondecreasing, resp.\ nonincreasing with
respect to this partial order.

\begin{prop}[Harris-Fortuin-Kasteleyn-Ginibre inequality]
If $\mathcal{A}$ and $\mathcal{B}$ are both increasing
or both decreasing events,
\begin{equation}
\label{FoKaGi}
\Q^{\jota} (\mathcal{A} \cap \mathcal{B}) \; \geq\;
\Q^{\jota}(\mathcal{A})\,\Q^{\jota}(\mathcal{B})\;.
\end{equation}
\end{prop}

In the context of continuous-time percolation, an inequality of this
sort was first proven by Bezuidenhout and Grimmett
(1991\nocite{bezgri91}).  They did so by showing that the process is
the weak limit of discrete-time processes that satisfy the
corresponding inequality.  Their approach provides at the same time a
proof of the van den Berg-Kesten inequality, but has the inconvenience
of imposing an additional topological condition on events, namely that
their boundary ---in some suitable metric topology--- have measure
zero.  While such requirement is satisfied for the events of interest
to us, a different approach, based on the martingale convergence of
expectations of the relevant events, yields the result as stated
above.  A proof along this line is provided in the book by Meester and
Roy (1996)\nocite{meeroy96} (see their Theorem 2.2) for the
percolation of spheres of random radii on a continuous space (Poisson
Boolean model).  The proof is easily adaptable to our setting of
percolating cylinders.
\smallskip

The second inequality refers to increasing events happening in a
disjoint manner.  For brevity, we state it only for the type of events
needed in the sequel.  Let us consider a box $B = \Lambda \times I$,
with $\Lambda\subset {\Z}^d$ and $I$ a closed finite interval in
${\R}$, and space-time regions $B_1,\ldots,B_n \subset B$.  We denote
$\{ B_1 \rightarrow_{B} B_2 \} \circ \cdots \circ \{ B_{n-1}
\rightarrow_{B} B_n \}$ the event of having $n-1$ open paths,
respectively connecting in $B$ some point in $B_{i}$ to some point in
$B_{i+1}$, $1\le i\le n-1$, such that no two paths share the same
cylinder.

\begin{prop}[particular case of the van den Berg-Kesten inequality]
For  $B$, $B_1,\ldots,B_n$ as above,
\begin{equation}
\label{BeKe}
\Q^{\jota}\Bigl\{\{ B_1 \rightarrow_{B} B_2 \} \circ \cdots \circ \{ B_{n-1}
\rightarrow_{B} B_n \}\Bigr\} \;\le\;
\Q^{\jota}\Bigl\{ B_1 \rightarrow_{B} B_2 \Bigr\} \cdots \,
\Q^\jota\Bigl\{ B_{n-1} \rightarrow_{B} B_n \Bigr\}\;.
\end{equation}
\end{prop}

This inequality is a consequence of the more general inequality proven
by van den Berg (1996)\nocite{ber96} for increasing events of marked
Poisson processes.  The only subtlety is that van den Berg's result
requires the events to depend on Poisson clocks ringing within a
\emph{bounded} region of $\R^{d+1}$.  The connectivity events in
\reff{BeKe} do not seem to fit in this framework because they may be
determined also by cylinders born in an arbitrarily remote past.  To
obtain \reff{BeKe} we must, therefore, apply van den Berg's result to
the marked process in $B$ obtained by adding, independently, an
invariant initial distribution of cylinders.  This corresponds to a
further, independent, spatial Poisson marked process on
$\Lambda\times\{\inf I\}$ with rates
$\{w^\jota(\gamma):\gamma\in\G_\Lambda\}$.  Less direct proofs are
also possible either by adapting the Bezuidenhout and Grimmett (1991)
approach mentioned above, or the proof of Theorem 2.3 of Meester and
Roy (1996).
\smallskip

Our last inequality is obtained in Klein (1994)\nocite{kle94} as a
corollary of van den Berg's inequality.  In our oriented setting, we
can present a totally different argument.

\begin{prop}[Hammersley-Lieb-Simon inequality]
Let
$L,T\in\R_+$,
  $\delta$ as defined below \reff{r.delta} and $W\subset {\Z}^d \times
  {\R}$.  Then, for every $X=(x,t_X)\in W$ and $Y=(y,t_Y) \in W
  \setminus B_{L + \delta,T}(X)$, with $t_Y\le t_X$, we have
\begin{equation}
\label{HSL2}
G_W^{\jota}(X,Y) \;\leq\; \sum_{Z\, \in\, \partial B_{L,T}(X)\,\cap\, W}
G_{B_{L+ \delta, T}(X)}^{\jota}(X,Z)\; G_W^{\jota}(Z,Y_Z)
\end{equation}
with $Y_Z=(y,\min(t_Y,t_Z))$.
\end{prop}

\proof The enclosing character of $B_{L+\delta,T}(X)$ and the choice
of $\delta$, implies that every connection from $X$ to $Y$ includes a
connection within $B_{L+\delta,T}(X)$ to a point $Z \in \partial
B_{L,T}(X)\,\cap\, W$. It is here where the assumption of bounded animal sizes
is critically needed.  This point $Z$ is subsequently joined by an
open path to a final cylinder whose section contains $y$.  If $t_Y\le
t_Z$, this open path determines the event $Z\to_W Y$ (recall that only
connections that are backwards in time are considered).  If on the
contrary $t_Z\le t_Y$, then this final cylinder must have been born
before $t_Z$ and stood alive at least up to $t_Y$.  In particular, it
contains the point $(y,t_Z)$.  Both situations are summed up in the
inequality
\begin{equation}\label{cam:20}
G_W^{\jota}(X,Y)\;\le\; \sum_{Z\,\in\,\partial B_{L,T}(X)\,\cap\, W}
\Q^\jota\Bigl\{X \rightarrow_{B_{L + \delta, T}(X)} Z\,;\,
Z \rightarrow_W Y_Z \Bigr\}\;.
\end{equation}
(The sum in the right-hand side involves, in fact, a time integral.
Its justification requires, therefore, a limit of conectivities
involving time-discretizations of $\partial_V B_{L,T}(X)$.  These
connectivities are continuous functions, so the limit converges to the
integral.  We omit the details.)

The desired inequality \reff{HSL2} follows from \reff{cam:20} and the
fact that the events $X\rightarrow_{B_{L + \delta, T}(X)}Z$ and $Z
\rightarrow_W Y_Z$ are \emph{independent}.  Indeed, they are
determined, respectively, by which cylinders are alive for $t<t_Z$ and
which are alive for $t>t_Z$.  Such events are independent due
of the exponential character of cylinder lives.  Perhaps the simpler
way to see this is by resorting to an alternative construction of the
cylinder process, where for each animal $\gamma$ the birth and death
events are generated independently with respective exponential rates
$w^{\jota}(\gamma)$ and $1$.  Cylinders are born at a birth event and live up
to the next death event, neglecting intermediate events of the wrong
type.  With this construction, connectivity events before and after a
given time $t_Z$ are determined by different, hence independent, birth
and death events.  \square

\begin{coro}\label{c.HSL}
  Let
$L,T\in\R_+$, $\delta$ as
  defined below \reff{r.delta} and $W\subset {\Z}^d \times \R$.  Then,
  for every $X=(x,t_X)$ and $Y=(y,t_Y)$ in $W$ with $t_Y\le t_X$, we
  have that,
\begin{eqnarray}
\label{HSL2.M}
G_W^{\jota}(X,Y) &\leq& \sum_{Z_1\, \in\, \partial B_{L,T}(X)\,\cap\, W}
G_{B_{L+ \delta, T}(X)}^{\jota}(X,Z_1)
\sum_{Z_2\, \in\, \partial B_{L,T}(Z_1)\,\cap\, W}
G_{B_{L+ \delta, T}(Z_1)}^{\jota}(Z_1,Z_2)\;\cdots\nonumber\\
&&\qquad\qquad\cdots\;
\sum_{Z_N\, \in\, \partial B_{L,T}(Z_{N-1})\,\cap\, W}
G_{B_{L+ \delta, T}(Z_{N-1})}^{\jota}(Z_{N-1},Z_N)
\; G_W^{\jota}(Z_N,Y_{Z_N})\nonumber\\
\
\end{eqnarray}
for every
\begin{equation}\label{HSL.N}
N \;\le\; \hbox{\rm integer part of }
\max \left( \frac{\norm{x-y}}{L+\delta}\,,\,
\frac{| t_X-t_Y |}{T} \right)\;.
\end{equation}

\end{coro}

\proof This is just an iteration of \reff{HSL2}.  The number of times
such iteration can be performed is at least equal to the right-hand
side in \reff{HSL.N}.
\square

\subsection{Regularity and ``good enough'' probabilities}\label{s.multi}

In this section we introduce the main notions defining the multiscale
approach.

From now on we consider boxes as in \reff{boxes} where the temporal
height $T$ is an increasing function of $L$.  In Section
\ref{s.change} the function $T(L)$ will be eventually chosen as a
stretched exponential, but the following results do not depend on such
a particular choice.  To simplify the notation we characterize the
boxes by its spacial lenght $L$ and denote
\begin{equation}
\label{sbox}
B_{L}(X) = B_{L+\delta, T(L)}(X).
\end{equation}

Following Klein (1994)\nocite{kle94} we separate ${\Z}^d$ in
regular and singular regions for a fixed realization of the
environment $\jota$.
\begin{defin}
Let $m>0$ and $L>1$. A site $x \in {\Z}^d$ is said to be
{\bf $(m,L)-$regular} if
\begin{equation}
\label{defreg}
G_{B_L(x,0)}^{\jota}\Bigl((x,0)\,,\, \partial\Bigr)
\;\leq\; e^{-m(L+\delta)}.
\end{equation}
Otherwise $x$ is called {\bf $(m,L)-$singular}. A set $\Lambda \subset
{\Z}^d$ is said $(m,L)-$regular if every $x \in \Lambda$ is
$(m,L)-$regular; otherwise it is $(m,L)-$singular.
\end{defin}

Regularity will be used in conjunction with the Hammersley-Lieb-Simon
inequality through the following crucial result.
\begin{lema}
\label{bueno}
Let $\Lambda$ be a $(m,L)$-regular region, $W \supset \Lambda \times
\R$.  Then, for every $X=(x,t_X)$ and $Y=(y,t_Y)$ in $W$,
\begin{equation}
\label{ragalar}
G_{W}^{\jota}(X,Y) \;\leq\; \exp\Bigl\{-m\,(L+ \delta)\,
N_\Lambda(X,Y)\Bigr\}
\end{equation}
with
\begin{equation}\label{n.ragalar}
N_\Lambda(X,Y) \;=\; \hbox{\rm integer part of }
\min\left\{\frac{{\rm dist}(x,\Lambda^c)}{L + \delta}\;,\;
 \max \left[ \frac{\norm{x-y}}{L+\delta}\,,\,
\frac{|t_X-t_Y|}{T} \right] \right\}.
\end{equation}
\end{lema}
\proof This follows from Corollary \ref{c.HSL}.  The value of
$N_\Lambda(X,Y)$ satisfies the constraint \reff{HSL.N} and guarantees
that all the intermediate sites $z_i$ are in $\Lambda$ and, hence, they are
regular. We then bound the right-hand side of \reff{HSL2.M} starting
from the right:  $G_W^{\jota}(Z_N,Y_{Z_N})$ is bounded by one, and each
of the preceding sums by $\exp\{-m\,(L+\delta)\}$.  \square
\medskip

We also formalize the notion of \emph{scale}
\begin{defin}
A {\bf scaled sequence} is a triple $(L_0,\alpha,T)$ where $L_0,
  \alpha>1$ and $T:(0,\infty)\to (0,\infty)$ is a function that grows
  faster than any power.  Each such triple defines an increasing
  sequence of sizes $L_{k+1}=L_{k}^{\alpha}$ for $k=0,1,2, \ldots$.  The
  length $L_k$ is the $k$-th scale of the sequence.
\end{defin}

Finally, we associate ``good probabilities'' to scales.
\begin{defin}
A scaled sequence $(L_0,\alpha,T)$ has {\bf $m_\infty$-good-enough
probabilities} if there exists $p>\alpha d$ such that
\begin{equation}
\label{pipio}
\P \Bigl\{x \mbox{ is } (m_{\infty},L_k)-\mbox{regular } \Bigr\}
\;\geq\; 1 - \frac{1}{L_k^p}
\end{equation}
for all $k=0,1,2, \ldots$, for all $x\in\Z^d$.
\end{defin}

We end this section with the proof of the ``easy part'' of the
multiscale argument.
\begin{teo}
\label{behavior}
If a scaled sequence $(L_0,\alpha,T)$ has $m_\infty$-good-enough
probabilities, then, for any $m\in(0,m_\infty)$ there exist constants
$\{C_x^\jota(m): x\in\Z^d\}$ with $\P\{C_x^\jota(m)<\infty:
x\in\Z^d\}=1$, such that
\begin{equation}\label{cam.30}
G^{\jota}(X,Y) \;\leq\; C_x^\jota(m)\,
 \exp\Bigl\{-m \max \Bigl[\norm{x-y}\,,\,
 T^{-1}(t_X-t_Y)\Bigr] \Bigr\}
\end{equation}
for all $X,Y \in {\Z}^d\times\R$ with $t_X\ge t_Y \in {\R}$.
\end{teo}

\proof We follow Klein (1994)\nocite{kle94}. Take $b>1$ (to be
determined later) and consider
\begin{equation}
{\cal S}_k \;:=\; \Bigl\{\Lambda[x;b(L_{k+1} + \delta)] \mbox{ is not a }
(m_{\infty},L_k)- \mbox{regular region} \Bigr\}.
\end{equation}
This event is verified if at least one of the sites in
$\Lambda[x;b(L_{k+1} + \delta)]$ is $(m_\infty,L_k)$-singular.  As, by
hypothesis, the probability for a given site to be
singular
is at most $L_k^{-p}$, we obtain
\begin{eqnarray*}
\PP\{ {\cal S}_k \} &\leq& \sum_{y \in \Lambda[x;b(L_{k+1}+ \delta)]} \,
\PP\Bigl\{ y \mbox{ is }(m_{\infty},L_k)- \mbox{singular }\Bigr\} \\[5pt]
&\leq& \frac{[2b(L_{k+1} + \delta)]^d}{L_k^p}
\;\leq\; \frac{(4bL_k^\alpha)^d}{L_k^p}
\;=\; \frac{(4b)^d}{L_0^{\alpha^k(p- \alpha d)}}\;.
\end{eqnarray*}
As $p> \alpha d$, this bound shows that the probabilities $\PP\{{\cal
  S}_k\}$ are summable in $k$.  Therefore, by Borel-Cantelli, with
probability $1$ there exists $k_1 =k_1(x,\ell_0, b,\jota) < \infty$,
such that $\Lambda[x;b(L_{k+1} + \delta)]$ is $(m_{\infty}, L_k)$-
regular region for all $k \geq k_1$.

Now fix $X$, and classify the sites $Y$ into regions
\begin{equation}\label{parak.1}
\overline{\mathcal{R}} \;=\; \Bigl\{ Y :\;
\max\Bigl[ \norm{x-y}\,,\,
T^{-1}(t_X-t_Y) \Bigr]\; <\; b(L_{k_1}+ \delta) \Bigr\}
\end{equation}
and
\begin{equation}
\label{parak}
\mathcal{R}_k \;=\; \Bigl\{ Y :\;
b(L_k + \delta) \;\leq\;  \max\Bigl[ \norm{x-y}\,,\,
T^{-1}(t_X-t_Y) \Bigr] \;<\; b(L_{k+1}+\delta)\Bigr\}
\end{equation}
for $k\ge k_1$.

Let us first consider $Y\in\mathcal{R}_k$ for some $k\ge k_1$.  In
this case we have $\norm{x-y} < b(L_{k+1}+\delta)$ and so
$y \in \Lambda[x;b(L_{k+1}+\delta)]$ which is a
$(m_{\infty},L_k)$-regular region. It follows, from Lemma \ref{bueno}
(for $W=\Z^d\times \R$ and $\Lambda=\Lambda[x;b(L_{k+1}+\delta)]$),
that
\begin{equation}\label{cam.31}
G^{\jota}(X,Y) \leq e^{-m_{\infty}(L_k + \delta)\,N}
\end{equation}
with
\begin{equation}
\label{enne}
N\;=\; \hbox{\rm integer part of }
\min \left\{ \frac{b(L_{k+1}+ \delta)}{L_k + \delta}, \max
\left\{ \frac{\norm{x-y}}{L_k + \delta}, \frac{t_X-t_Y}
{T(L_k)}\right\} \right\}.
\end{equation}

A first bound of \reff{cam.31} comes from the observation that, as
$Y\in\mathcal{R}_k$ implies
$\norm{x-y}\le b(L_{k+1}+\delta)$,
\[
N \;\ge\; \hbox{\rm integer part of }
\frac{\norm{x-y}}{L_k + \delta}\;.
\]
Therefore,
\begin{equation}
\label{espacio}
G^{\jota}(X,Y) \;\leq\; \exp\left\{-m_{\infty}(L_k + \delta) \left(
    \frac{\norm{x-y}}{L_k + \delta} -1 \right)\right\}
\end{equation}
(we bounded the integer part of a number by the number minus one).

This bound can be improved in cases where the temporal part dominates
in the sense that $T^{-1}(t_X-t_Y) >\, \norm{x-y}$.
In this case we use that $b(L_k + \delta) \leq T^{-1}(t_X-t_Y)$
for $Y\in\mathcal{R}_k$, so that, as $T$ grows faster than any power,
\[
\frac{t_X-t_Y}{T(L_k)} \;\geq\; \frac{T(b(L_k + \delta))}{T(L_k)}
\;\geq \; \frac{b(L_{k+1}  + \delta)}{L_k + \delta}
\;\geq\; \frac{\norm{x-y}}{L_k + \delta}
\]
if $k_1$, and hence $k$, is chosen large enough.  Hence
\begin{equation}\label{cam.35}
N \;\ge\; \hbox{\rm integer part of }
\frac{b(L_{k+1}+ \delta)}{L_k + \delta}
\;\ge\; \hbox{\rm integer part of }
\frac{T^{-1}(t_X-t_Y)}{L_k + \delta}
\end{equation}
where the last inequality comes from the definition of
$\mathcal{R}_k$.  From this and \reff{cam.31} we obtain that
\begin{equation}
\label{tiempo}
G^{\jota}(X,Y) \;\leq\; \exp\left\{-m_{\infty}(L_k + \delta)\left(
    \frac{T^{-1}(t_X-t_Y)}{L_k + \delta} -1 \right)\right\}
\end{equation}
whenever $T^{-1}(t_X-t_Y) >\, \norm{x-y}$.

Inequalities (\ref{espacio}) and (\ref{tiempo}) can be combined in the
expresion
\begin{equation}\label{cam.36}
G^{\jota}(X,Y) \leq \exp \Bigl\{ -m_{\infty} \Bigl[ \max\Bigl(
\norm{x-y}\,,\, T^{-1}(t_X-t_Y) \Bigr)
 - (L_k + \delta) \Bigr] \Bigr\}\;.
\end{equation}
The additive correction $(L_k+\delta)$ can be turned into a factor
$(1-b^{-1})$ because $(L_k + \delta) \leq (1/b) \max\{ \norm{x-y},
T^{-1}(t_X-t_Y) \}$, by definition of $\mathcal{R}_k$.  In addition we
choose
$b = \frac{m_{\infty}}{m_{\infty} -m}$
so that $m = m_{\infty}(1- \frac{1}{b})$.  In this way \reff{cam.36}
yield
\begin{equation}\label{cam.37}
G^{\jota}(X,Y) \;\leq\; \exp\Bigl\{-m \max\Bigl[ \norm{x-y}
\,, \,  T^{-1}(t_X-t_Y) \Bigr]\Bigr\}
\end{equation}
uniformly in $k\ge k_1$.

Finally we consider $Y\in\overline{\mathcal{R}}$.  They satisfy
$\max\{ \norm{x-y}, T^{-1}(t_X-t_Y) \} < b\,(L_{k_1}+
\delta)$, so we can write
\begin{eqnarray}
\label{Aquiles}
G^{\jota}(X,Y) &\leq&  G^{\jota}(X,Y)\,
\exp\Bigl\{mb\,(L_{k_1} + \delta )\Bigr\}
\, \exp\Bigl\{- m \max\Bigl[ \norm{x-y}\,,\,
T^{-1}(t_X-t_Y) \Bigr]\Bigr\}
\nonumber\\
&\leq& C_x^\jota(\ell_0,m)\, \exp\Bigl\{- m \max\Bigl[ \norm{x-y},
T^{-1}(t_X-t_Y) \Bigr]\Bigr\}
\end{eqnarray}
where $C_x^\jota(\ell_0,m) = e^{mb(L_{k_1}+ \delta)}$.  The desired
bound \reff{cam.30} follows from \reff{cam.37} and \reff{Aquiles}.
\square

\subsection{The change of scale}\label{s.change}

\subsubsection{The change-of-scale theorem}
This section contains the heart of the multiscale argument.  Its main
result is the following theorem, related to Theorem 3.2 of Klein
(1994)\nocite{kle94}, which establishes sufficient conditions for a
scaled sequence to have good-enough probabilities.

\begin{teo}
  \label{porfin}
  Let $d\geq1$ and choose $a$ satisfying \reff{eq:r2}.  Assume that
  the disorder is such that
 $\aleph<\infty$ satisfies \reff{eq:rj1}
  for the chosen $a$.  Let $\alpha = d + \sqrt{d^2 + d}$ and take
  $\nu$ and $p$ such that
\begin{equation}
\label{pf3}
\frac{\alpha d\,[\alpha + a +1]}{a\,[\alpha - d + \alpha d]} \;< \nu \;
\end{equation}
and
\begin{equation}
\label{pf4}
\alpha\, d \;<\; p \;<\;
\frac{a(\nu[\alpha - d + \alpha d] - \alpha d) - \alpha d}{\alpha} \;.
\end{equation}
Finally, consider $m_0$ and $m_{\infty}$ such that $0< m_{\infty} <
m_0$. Then, there exists
$\widetilde{L} = \widetilde{L}(d,\delta,a,
\aleph, \nu,p,m_0,m_{\infty})< \infty$ such that if for some
$L_0>\widetilde L$
\begin{equation}
\label{pf5}
\P\Bigl\{ x \mbox{ is } (m_0 ,L_0)-\mbox{regular} \Bigr\}
\;\geq\; 1 - \frac{1}{L_0^p}\;,
\end{equation}
for all $x\in\Z^d$, then the scaled sequence $(L_0,\alpha,e^{L^\nu})$ has
$m_\infty$-good-enough probabilities.
\end{teo}

Note that condition \reff{pf1} guarantees the existence of $\nu$
satisfying \reff{pf3}, which in turns implies the existence of $p$ as
in \reff{pf4}.

The combination of this theorem and Theorem \ref{behavior}, implies
that the key Lemma \ref{main} ---and hence all the properties stated
in Theorem \ref{225}--- follow once \reff{pf5} is satisfied.  Note
that the logarithmic time-dependence in the key lemma is a consequence
of the choice
\begin{equation}\label{oud.1}
T(L) \;=\; \exp(L^\nu)\;,
\end{equation}
$0<\nu<1$, for the time-height of the scaled sequence.

The proof of Theorem \ref{porfin} is inductive.  We show
that it is possible to arrange the parameters in such a way that, if
\begin{equation}
\label{induc1}
\P\Bigl\{ x \mbox{ is } (m,l)- \mbox{regular } \Bigr\}
\;\geq\; 1 - \frac{1}{l^p}
\end{equation}
for $x \in \Z^d$, $m>0$ and $l$ sufficiently large, then we also have
\begin{equation}
\label{induc2}
\P\Bigl\{ x \mbox{ is } (M,L)- \mbox{regular } \Bigr\} \;\geq\;
1 - \frac{1}{L^p}
\end{equation}
where $L = l^{\alpha}$ and $M$ is an appropiate function of $m$.  The
inductive step has two ingredients.  First, in Subsection \ref{s.geo}
we show that the inductive regularity holds in the presence of two
events $\mathcal A$ and $\mathcal{B}_\Delta$.  This part of the
argument is based on geometric considerations and the inequalities of
Section \ref{s.tool}.  Subsequently, in Subsection \ref{s.pro} these
events are proven to hold with high-enough probability.  The success
of the approach relies on the careful definition of the events
$\mathcal A$ and $\mathcal{B}_\Delta$.

\subsubsection{Geometrical estimates:  Good events imply good
  behavior}\label{s.geo}

\begin{lema}
\label{geoest}
Consider some fixed $x\in\Z^d$.
Let $\nu,\alpha,l$ be such that $0<\nu<1$, $\alpha>1$ and
$l>\delta$. Set
\begin{equation}
  \label{eq:rr.30}
\theta_0 \;:=\; \min\{\alpha -1, \alpha(1-\nu)\}
\end{equation}
and take $m_0,\theta, m$ with $m_0>0$, $0<\theta<\theta_0$
and
$l^{-\theta}<m<m_0$. Put $L=l^{\alpha}$, pick a positive integer $R$
and define the event
\begin{eqnarray}
\label{star1}
{\cal A} &=& \Bigl\{\jota : \mbox{there exists } x_1, x_2, \ldots ,
x_R \in \Lambda[x; L + \delta] \mbox{ such that } \nonumber \\
& & \Lambda[x;L + \delta] \setminus \bigcup_{j=1}^{R} \Lambda[x_j ;
2(l + \delta)+1] \mbox{ is a } (m,l)-\mbox{regular region}
\Bigr\}\;.
\end{eqnarray}
Assuming ${\cal A}$ is true, take $\Delta>0$, $b$ such that $0<b< \alpha
\nu$ and $\kappa > \max \{1, \nu + \theta_0 \}$.  Define the event
\begin{equation}
  \label{eq:r20} {\cal B}_{\Delta} \;=\; \Bigl\{ \jota: \prod_{\gamma
\in \ani_{\widetilde{\Lambda}}}
K_{\Delta}^{\jota}(\gamma) \,\ge\, e^{-l^b} \Bigr\}\;,
\end{equation}
where
\begin{equation}
\label{tildado}
 \widetilde{\Lambda} \;=\; \bigcup_{j=1}^{R} \Bigl(\Lambda[x_j;l^{\kappa}]
\bigcap \Lambda[x;L + \delta] \Bigr)
\end{equation}
---the union being taken on the points involved in the definition of the
  event ${\cal A}$--- and
\begin{equation}
\label{hache}
K_{\Delta}^{\jota}(\gamma) \;=\;  e^{-(1+\Delta )w^{\jota}(\gamma)} +
(1- e^{-w^{\jota}(\gamma)})(1 - e^{- \Delta}) \, e^{-\Delta
  w^{\jota}(\gamma)}\;.
\end{equation}
Then, there exist $a_0 = a_0(d, \delta, \alpha,\nu,m_0,R)\ge 0$ and $l_0 =
l_0(d,\delta,\alpha,\nu,\kappa,m_0,\theta,\theta_0,R) < \infty$ such
that if $l>l_0$ the following holds: If ${\cal A}$ and ${\cal
  B}_{\Delta}$ are true for some $\Delta \in (0,1]$ then the
site
$x \in \Z^d$ is $(M,L)$-- regular with
\begin{equation}\label{rr:5}
M\;\geq\; m - \frac{a_o}{l^{\theta_0}} \;\geq\; \frac{1}{L^\theta}\;.
\end{equation}
\end{lema}

This lemma is a direct consequence of the following Sublemas
\ref{sl.first} and \ref{sl.second}, where we analyze separately the
connectivity of the site
$X=(x,0)$ to sites in the vertical and horizontal
boundary.  The regularity stated below \reff{hache} follows
immediately by summing the corresponding bounds over both boundaries.

Sublemma \ref{sl.first} involves estimations purely in the spatial
direction, hence only the event $\mathcal{A}$ is relevant and the
choice of the time scale $T(L)$ actually plays no role
(though, for concreteness, we stick to the subexponential dependence).
Event $\mathcal{B}_\Delta$ is required to control the time-like
percolation studied in Sublemma \ref{sl.second}.  It provides a lower
bound for the probability of not having long towers of ancestors based
around ``deffective'' sites [see \reff{eq:rr.23}--\reff{eq:rr.22}
below].
The stretched-exponential choice \reff{oud.1} becomes essential for
this sublemma which is in fact the hardest estimation of the paper.

\begin{sublema}\label{sl.first}
There exist $l_1 = l_1(d,\delta,\alpha,\nu,m_0,\theta,R) < \infty$
and
$a_1 = a_1(d, \delta, \alpha,\nu,m_0,R)>0$ such that for all $l
  > l_1$, all $Y \in \partial_V B_L(X)$ and all $\jota\in{\cal A}$
\begin{equation}\label{rr:3}
 G_{B_L(X)}^{\jota}(X,Y) \;\leq\; e^{-M_1 (L+\delta)}
\end{equation}
with
\begin{equation}
M_1 \;=\; m - \frac{a_1}{l^{\theta_0}} \;\geq\;
\frac{1}{L^{\theta}}\;.
\end{equation}
\end{sublema}

\proof We work with a fixed environment $\jota \in {\cal A}$.  First we
group the deffect-sets $\Lambda[x_j;2(l+\delta) +1]$ into larger cubes
which absorb ``pockets'' totally surrounded by original deffects.
This leaves us with a much simpler geometrical situation, where the
regular region is the complement of a finite family of cubes.  More
precisely, elementary geometrical considerations show that there
exists a possibly smaller collection of sites $y_1, y_2, \ldots,
y_{R'} \in \Lambda[x;L+ \delta]$, with $R' \leq R$, $n_1, n_2, \ldots
, n_{R'} \in \{1,2, \ldots, R\}$ and $n_1 + n_2 + \dots + n_{R'} \leq
R$, such that the sets $\Lambda[y_i; n_i(2(l+ \delta) +1)]$, $i=1,2,
\ldots, R'$, are at a distance strictly larger than one,
$\bigcup_{j=1}^{R} \Lambda[x_j;2(l +\delta) +1] \subset
\bigcup_{i=1}^{R'} \Lambda[y_i;n_i(2(l + \delta)+1)]$, and $\Lambda'
:= \Lambda[x;L+\delta] \setminus \bigcup_{i=1}^{R'} \Lambda[y_i;
n_i(2(l +\delta)+1) ]$ is a \emph{connected} $(m,l)$-regular region.

These problematic regions define cylinders of radius $n_i(2(l+
\delta)+1)$ and temporal height $e^{L^{\nu}}$ centered at the
points
$y_i$:
\[
B_i \;=\; B_{n_i(2(l+\delta)+1),e^{L^{\nu}}}((y_i,0))\;, \quad
i=1,2, \ldots , R'
\]
[recall the notation in \reff{boxes}].
We shall control the connectivity function in the spacial direction,
on the $(m,l)$-regular region
\begin{equation}
B'\;=\; B_L(X) \setminus  \bigcup_{i=1}^{R'} B_i
\end{equation}
[we use the notation resulting from \reff{sbox} and \reff{oud.1}].

If $X \in B'$, we denote $\partial B_0 =\{X\}$; otherwise $X \in
B_{i'}$ for some $i'$ and we put $\partial B_0= \partial
B_{i'}$. Similarly, if $Y \in B'$, we denote $\partial B_{R+1}
=\{Y\}$; otherwise $Y \in B_{i''}$ for some $i''$ and we put $\partial
B_{R+1}= \partial B_{i''}$.  Every connection from $X$ to $Y$ can be
decomposed into disjoint connections among some of the cylinders
$B_i$.  Therefore,
\begin{eqnarray}
\label{patroclo}
\Bigl\{ X \rightarrow_{B_L(X)} Y \Bigr\} &\subseteq& \bigcup_{r=1}^{R'}\;
 \bigcup_{\{i_1, i_2,  \ldots , i_r\} \subset \{1,2, \ldots ,R'\}}
\Bigl\{ \partial B_0 \rightarrow_{B'}  \partial B_{i_1} \Bigr\}
 \circ \Bigl\{ \partial B_{i_1} \rightarrow_{B'} \partial B_{i_2}
 \Bigr\} \nonumber \\
&& \qquad\qquad\qquad \qquad\qquad\qquad \circ
  \cdots  \circ  \Bigl\{ \partial B_{i_r} \rightarrow_{B'}
\partial B_{R'+1}\Bigr\}\;.
\end{eqnarray}

As $B' \subset \Lambda' \times \R$ with $\Lambda'$ a
$(m,l)$-regular region, we can apply \reff{ragalar} to obtain
\begin{equation}\label{rr:1}
{\Q}^{\jota}\{ \partial B_{j_1} \rightarrow_{B'} \partial B_{j_2} \}\;
\leq\;
\Bigl[(R[2(l + \delta) +1])^d e^{L^{\nu}} \Bigr]^2 \exp\Big\{ -m(l+
\delta) \Big( \frac{D_{j_1,j_2}}{l+\delta} -1 \Big)\Big\}\;.
\end{equation}
We have denoted
\[
D_{j_1,j_2} \;=\; \min \Bigl\{ \norm{x_1 - x_2}\;: \;
(x_1,t) \in \partial B_{j_1}\,, \,(x_2,s) \in
\partial B_{j_2} \mbox{ for some } t,s \Bigr\}
\]
and used the bounds
\[
\max\left\{ \frac{\norm{x_{j_1} - x_{j_2}}}{l+\delta}\,,
  \frac{|t-s|}{e^{L^{\nu}}} \right\} \;\geq\;
 \frac{\norm{x_{j_1} - x_{j_2}}} {l+\delta} \;\geq\;
  \frac{D_{j_1,j_2}}{l+\delta}\;.
\]

>From \reff{patroclo}, \reff{rr:1} and the van den Berg-Kesten
inequality \reff{BeKe} we have
\begin{eqnarray}
G_{B_L(X)}^{\jota} (X,Y) &\leq& \sum_{r=1}^{R'} \sum_{\{i_1, i_2, \ldots , i_r\}
  \subset  \{1,2, \ldots ,R'\}} \,
\Bigl[(R[2(l + \delta) +1])^d e^{L^{\nu}} \Bigr]^{2(R'+1)} \nonumber \\
& & \times \exp\left\{ -m(l+ \delta)\left( \frac{D_{0,i_1}+D_{i_1,i_2}+ \ldots
  + D_{i_r, R'+1}}{l+\delta} -(r+1) \right) \right\}. \nonumber
\end{eqnarray}
This inequality, together with the bounds $D_{0,i_1}+D_{i_1,i_2}+
\ldots + D_{i_r,R'+1} \geq L - 2[2(l+\delta)+1] [n_{i_1} + n_{i_2} +
\ldots + n_{i_r}] \geq L - 2[2(l+\delta)+1]R$ and $ \sum_{r=1}^{R'}
\sum_{\{i_1, i_2, \ldots , i_r\} \subset \{1,2, \ldots ,R'\}} \, 1
\leq (R+1)!$, yield
\begin{eqnarray}\label{rr:2}
G_{B_L(X)}^{\jota}(X,Y) &\leq& (R+1)!
\Bigl[(R[2(l + \delta) +1])^d e^{L^{\nu}} \Bigr]^{2(R+1)} \nonumber \\
& & \times  \exp\left\{
  - m(l+ \delta) \left( \frac{L- 2[2(l+\delta)+1]R}{l+\delta} - (R+1) \right)
\right\}\;.
\end{eqnarray}
Therefore
\[ G_{B_L(X)}^{\jota}(X,Y) \leq e^{-M'_1 (L + \delta)};\]
with
\[ M'_1 \;=\;
m - \frac{c_1m_0}{l^{\alpha -1}} - \frac{c_2}{l^{\alpha(1-\nu)}} \]
for some fixed constants $c_1,c_2>0$ depending on $d$, $\delta$,
$\alpha$, $\nu$ and $R$.  This yields the proposed bound \reff{rr:3}
because, as $m \geq \frac{1}{l^{\theta}}$,
\[ M'_1 \;\geq\; M_1= m - \frac{a_1}{l^{\theta_0}}
\;\geq\; \frac{1}{L^{\theta}} \]
for some $a_1 = a_1(d,\delta,\alpha,\nu,m_0,R) >0$ and for $l$ large
enough (the meaning of ``large enough'' depends on
$d,\delta,\alpha,\nu,m_0,\theta$ and $R$). $\Box$
\medskip

\begin{sublema}\label{sl.second}
  Suppose the events ${\cal A}$ and ${\cal B}_{\Delta}$ in
  (\ref{star1}) and (\ref{eq:r20}) are true. Pick $\tau$ such that
  $\nu < \tau < \min\{\kappa - \theta_0 , \alpha \nu\}$. Then there
  exists $l_2 = l_2(d, \delta, \alpha,\nu,m_0, \theta,\kappa,
  b,\tau,R) < \infty$, such that for $l>l_2$ we have
\[ G_{B_L(X)}^{\jota}(X,Y) \leq exp\{ -M_2 e^{\frac{l^{\tau}}{4}} \} \]
for all $Y \in \partial_H B_L(X)$, with
\begin{equation}\label{rr:8}
M_2 \;=\; m - e^{-\frac{l^{\tau}}{3}} \;\geq\;
\frac{1}{L^{\theta}}\;.
\end{equation}
\end{sublema}

Note that it is possible to choose $\tau$ satisfying the hypothesis
because $\kappa > \nu + \theta_0$ and $\nu < \alpha \nu$.
\medskip

\proof
Let us fix an environment $\jota \in {\cal A}\cap {\cal
  B}_{\Delta}$.  The proof relies on the introduction of additional
space and time scales:
\begin{equation}
  \label{eq:rr.10}
  \widehat \delta (l) = l^\kappa \quad\mbox{and}\quad
\widehat\tau(l) = \exp(l^\tau)\;.
\end{equation}
The time scale $\widehat\tau(l)$ is intermediate between $T(l)$ and
$T(L)=T(l^\alpha)$ because $1<\tau<\alpha\nu$.  The space scale
$\widehat\delta$ is used to construct the region $\widetilde\Lambda$
defined in \reff{tildado}.  For the time being it is only required to
grow faster than $l$ ($\kappa>1$), but the hypotheses of Lemma
\ref{varsovia2} below will force $\kappa$ to be strictly smaller than
$\alpha$ [see condition \reff{condicion1}], making $\widehat\delta$ a
scale intermediate between $l$ and $L=l^\alpha$.
The intermediate time scale defines a partition of the box $B_L(X)$
into
\begin{equation}
  \label{eq:rr.11}
N \;=\; \hbox{integer part of }\; \frac{T(L)}{\widehat{\tau}(l)}
\end{equation}
``slices'' of time-heigth $\widehat{\tau}(l)$:
\begin{equation}
  \label{eq:rr.12}
S_j \;=\; B_{L+ \delta ,\widehat{\tau}(l)}\Bigl((x,-(j -
1)\widehat{\tau}(l)) \Bigr) \;, \quad j=1,2,
\ldots , N\;.
\end{equation}
We note that events in different $S_j$ are independent (see the
paragraph preceding Corollary \ref{c.HSL}).

We introduce the notation $B_{\Lambda} = \Lambda \times [-T(L),0]$ for
any $\Lambda \subset \Z^d$.  Let
\begin{equation}
  \label{eq:rr.16}
\widehat{\Lambda} \;=\; \Big( \bigcup_{j=1}^R \Lambda[x_j;2(l+ \delta)
+1] \Big) \, \bigcap\, \Lambda[x;L+\delta]
\end{equation}
be the set of ``deffective'' or ``irregular'' sites.  The region
$\widetilde\Lambda$ defined in \reff{tildado} corresponds to a
``buffer zone'' around, and including, these irregular sites.  We
decompose the event $\{ X \rightarrow_{B_L(X)} Y \}$ according to
whether
there exists at least one slice such that no deffective point
within it is connected to a preceding slice.  That is, we write
\begin{equation}
  \label{eq:rr.15}
 \Bigl \{ X \rightarrow_{B_L(X)} Y \Bigr\}\;=\;
\Bigl(\Bigl\{ X \rightarrow_{B_L(X)} Y \Bigr\}\cap A\Bigr) \,\bigcup\,
\Bigl(\Bigl\{ X \rightarrow_{B_L(X)} Y \Bigr\}\cap A^c \Bigr)
\end{equation}
with
\begin{equation}
  \label{eq:rr.16.1}
  A\;=\;\bigcup_{j=1}^N A_j \;=\;
\bigcup_{j=1}^N \Bigl\{B_{\widehat\Lambda}\cap S_j \longrightarrow_{S_j}
S_{j+1} \Bigr\}^c\;.
\end{equation}

The probability of the first event on the right of \reff{eq:rr.15} is
damped by regularity.  Indeed, the occurrence of $ \{ X
\rightarrow_{B_L(X)} Y \} \cap A $ implies the existence of a
connection in $B_{\Lambda[x;L+ \delta] \setminus \widehat{\Lambda}}$
of vertical height at least equal to $\widehat\tau(l)$.  Therefore,
\begin{equation}
  \label{eq:rr.17}
\{ X \rightarrow_{B_L(X)} Y \} \cap A  \;\subset\;
\bigcup \Bigl\{(y_1,s_1)
\longrightarrow_{B_{\Lambda[x;L+\delta]\setminus
  \widehat {\Lambda}}} (y_2,s_2) \Bigr\}\;,
\end{equation}
the union being taken over all $(y_1,s_1), (y_2,s_2) \in
B_{\Lambda[x;L+\delta] \setminus \widehat{\Lambda}}$ with $\mid s_1 -
s_2 \mid \geq \widehat\tau(l)$. By regularity (Lemma \ref{bueno})
we get, after simple computations
\begin{eqnarray}
\label{seco2}
\lefteqn{
{\Q}^{\jota} \Bigl(\{ X \rightarrow_{B_L(X)} Y \} \cap A\Bigr)
}\nonumber\\
&&  \leq\; \exp \Bigl\{ -m(l + \delta) \Bigl[
\frac{\widehat{\tau}(l)}{T(l)} -1 \Bigr] + 2d \ln(2L + 2
\delta)  + 2 \ln T(L)\Bigr\}\;.
\end{eqnarray}
Therefore,
\begin{eqnarray}
\label{secazo1}
\lefteqn{
\hspace{-1cm}
{\Q}^{\jota} \Bigl(\{ X \rightarrow_{B_L(X)} Y \} \cap A\Bigr)
}\nonumber\\
&\leq& \exp \Bigl\{ -m(l + \delta) \Bigl[
e^{l^{\tau} - l^{\nu}} -1 \Bigr] + 2d \ln(2l^{\alpha} + 2
\delta)  + 2 l^{\alpha \nu} \Bigr\} \nonumber \\
&\leq& \exp \Bigl\{- \bigl(m -e^{- \frac{l^{\tau}}{3}} \bigr)
\,e^{\frac{l^{\tau}}{2}} \Bigr\},
\end{eqnarray}
for $l$ large enough (recall that $\nu < \tau$).

To prove that under our hypotheses $A^c$ is sufficiently improbable,
we bound
\begin{equation}
  \label{eq:rr.18}
A_j \;\supset\; F_j \cap D_{-(j-1/2)\widehat{\tau}(l)}\;,
\end{equation}
where $F_j$ is the event that there is no connection inside the slice
$S_j$ from points in the deffective set to regular points outside the
buffer zone $\widetilde\Lambda$:
\begin{equation}
  \label{eq:rr.19}
F_j \;=\; \Bigl\{ \partial B_{\widehat{\Lambda}} \cap S_j
\longrightarrow_{S_j \setminus B_{\widehat{\Lambda}}}
B_{\Lambda[x;L + \delta] \setminus
  \widetilde{\Lambda}} \cap S_j  \Bigr\}^c\;,
\end{equation}
and $D_s$ represents the lack of connection within a strip of height
$\Delta$ involving deffective sites:
\begin{equation}
  \label{eq:rr.20}
  D_s \;=\; \left\{ \widetilde{\Lambda} \times \{ s \}
\longrightarrow_{\widetilde{\Lambda} \times [s- \Delta, s]}
  \widetilde{\Lambda } \times \{ s - \Delta \} \right\}^c\;.
\end{equation}
We are assumming that $l$ is sufficiently large
so that $\widehat\tau(l)/2>1$, and hence $\widehat\tau(l)/2>\Delta$
for all $\Delta\in(0,1]$.

Since each $F_j^c$ involves connections within regular regions, we can
apply Lemma \ref{bueno} to obtain
\begin{equation}
\label{qudefe}
{\Q}^{\jota} \{F_j^c \} \leq e^{-h(l)},
\end{equation}
with
\begin{equation}
\label{haga}
h(l) \;=\; m \{ \widehat\delta(l) -(3 l + 4 \delta +1) \} - 2 \ln
\widehat{\tau}(l) - d \ln\widehat\delta(l) - d \ln (2 l + 3 \delta +1)
- 2 \ln R\;.
\end{equation}
Hence
\begin{equation}
\label{hagaono}
h(l) \;\geq\; m\, c_1\, l^{\kappa} - c_2\, l^{\tau} \;\geq\;
 c_3\, l^{\kappa - \theta} \;\geq\; c_3\, l^{\kappa - \theta_0}\;,
\end{equation}
where $c_1,c_1$ and $c_3$ are positive constants ---depending on
$d,\delta,\kappa, \tau, \theta$ and $R$--- and $l$ is large enough. We
have used that $m \geq l^{-\theta}$, $\theta < \theta_0$ and $\tau <
\kappa - \theta$.

On the other hand,
\begin{equation}
  \label{eq:rr.21}
D_s \;\supset\; \bigcap_{\gamma \in \ani_{\widetilde{\Lambda}}} \,
\{ [ {\cal E}_1(\gamma)  \cup {\cal E}_2(\gamma)] \cap {\cal
  E}_3(\gamma) \}\;,
\end{equation}
where
\begin{eqnarray*}
{\cal E}_1(\gamma) &=& \{ \gamma \mbox{ is not present at time }
s-\Delta \}\;, \\
{\cal E}_2(\gamma) &=&  \{ \gamma \mbox{ is present at time } s-\Delta
\mbox{ but it does not survive until time } s \}\;, \\
{\cal E}_3(\gamma) &=& \{ \mbox{ there is no birth of } \gamma
\mbox{ in the interval } [s-\Delta , s] \}\;.
\end{eqnarray*}
We observe that
\begin{eqnarray}
\label{rr:7}
\Q^{\jota} \Bigl( [{\cal E}_1(\gamma) \cup {\cal E}_2(\gamma)] \cap {\cal
  E}_3(\gamma) \Bigr) &=& \Bigl(\Q^{\jota}[{\cal E}_1(\gamma)] +
  \Q^{\jota}[{\cal E}_2(\gamma)] \Bigr)\,
\Q^{\jota}[{\cal E}_3(\gamma)] \nonumber \\
&=& \Bigl[e^{-w^{\jota}(\gamma)} + (1- e^{-w^{\jota}(\gamma)})(1- e^{-
  \Delta})  \Bigr]
  e^{- \Delta w^{\jota}(\gamma)} \nonumber\\
&=& K_{\Delta}^{\jota}(\gamma)\;.
\end{eqnarray}
This result follows from well-known properties on Poisson processes
and from the fact that the events in which an animal $\gamma$ is
activated at time $-t$ and survives until time $0$ define an
inhomogeneous Poisson process of rate $w^{\jota}(\gamma)e^{-t}$.
Therefore, due to the $\Q^\jota$-independence of events pertaining
to different animals,
\begin{equation}
\label{eq:rr.23}
{\Q}^{\jota}\{D_s\} \;\ge\; \prod_{\gamma \in \ani_{\widetilde{\Lambda}}}
K_{\Delta}^{\jota}(\gamma) \;\geq\; e^{-l^b}\;.
\end{equation}

By the FKG-inequality we get, from \reff{eq:rr.18}, \reff{qudefe}
and
\reff{hagaono}--\reff{rr:7},
\begin{eqnarray}
  \label{eq:rr.22}
{\Q}^{\jota}\{A_j\} &\geq& {\Q}^{\jota}\{ F_j \}\;{\Q}^{\jota}\{
D_{-(j-1/2)\widehat{\tau}(l)} \}\nonumber\\
&\geq& \Bigl( 1 - e^{-c_3 l^{\kappa - \theta_0}} \Bigr)\, e^{-l^b}
\nonumber\\
& \geq& e^{-2l^b}
\end{eqnarray}
for $l$ large enough.
Thus, the independence of the events $A_j$, $j=1,2, \ldots ,N$ leads
to the bound
\begin{equation}
\label{acomp}
{\Q}^{\jota}(A^c) \;=\; \prod_{j=1}^N \Bigl(1- {\Q}^{\jota}(A_j)\Bigr)
\;\leq\; \bigl(1-e^{-2l^{b}} \bigr)^N
\;\leq\; \exp\Bigl(-N\,e^{-2l^b}\Bigr)\;.
\end{equation}
The condition $\tau < \alpha \nu$ implies that, for $l$ large enough,
$N \geq \exp(l^{\alpha \nu}/2)$, and as, besides, $b< \alpha \nu$,
we obtain from \reff{acomp} that
\begin{equation}
\label{acomple}
{\Q}^{\jota}(A^c) \;\leq\; \exp\Bigl(-e^{l^{\alpha \nu}/4}\Bigr)\;,
\end{equation}
for $l$ sufficiently large.

Finally, inequalities \reff{secazo1} and \reff{acomple} together with
the decompostion \reff{eq:rr.15} yield (recall that $\tau < \alpha
\nu$)
\begin{eqnarray}
G_{B_L(X)}^{\jota}(X,Y) &\leq&
{\Q}^{\jota} \Bigl(\{X \rightarrow_{B_L(X)} Y\} \cap A \Bigr )
+ {\Q}^{\jota} ( A^c ) \nonumber \\
&\leq& \exp \Bigl\{- \bigl(m -e^{- \frac{l^{\tau}}{3}} \bigr)
\,e^{\frac{l^{\tau}}{2}} \Bigr\} +
\exp \Bigl\{ - e ^{\frac{l^{\alpha \nu}}{4}}
\Bigr\} \nonumber\\
&\leq&  \exp \Bigl\{- \bigl(m -e^{- \frac{l^{\tau}}{3}} \bigr)
e^{\frac{l^{\tau}}{4}} \Bigr\}\;,
\end{eqnarray}
with $m -e^{- l^{\tau}/3} \geq 1/L^{\theta}$, for $l$ larger than a
certain threshold that depends on $d$, $\delta$, $\alpha$, $\nu$,
$m_0$, $\theta$, $\kappa$ ,$b$, $\tau$ and $R$. $\Box$

\subsubsection{Probabilistic estimates: Good events have high
  probability}\label{s.pro}

In this section we show that, assuming (\ref{induc1}), the events ${\cal A}$
and ${\cal B}_{\Delta}$ introduced in (\ref{star1}) and (\ref{eq:r20}) have
good enough probability for $l$ sufficiently large.  Hypothesis
\reff{eq:rj1} is used to bound the probability of the second event.
\begin{lema}
\label{varsovia1}
Let $m>0$, $l>\delta$, $\alpha >1$ and $p> \alpha d$. Put $L =
l^{\alpha}$. Pick a positive integer $R > \alpha p /(p - \alpha d)$
and define the event ${\cal A}$ as in (\ref{star1}). If
\begin{equation}
  \label{eq:rr.25}
\P \Bigl\{ x \mbox{ is } (m,l)- \mbox{regular } \Bigr\} \;\geq\;
  1 - \frac{1}{l^p}
\end{equation}
for all $x \in \Z^d$, then
\begin{equation}
\PP\{{\cal A}\} \;\geq\; 1 - \frac{1}{2L^p},
\end{equation}
assuming $l$ large exceed a certain minimum value that depends on $d$,
$\alpha$, $p$ and $R$.
\end{lema}

\proof Let us say that two sites $x_1,x_2 \in {\Z}^d$ are
$l$-nonoverlapping if the distance between the boxes
$\Lambda[x_1;l+\delta]$ and $\Lambda[x_2; l+\delta]$ is strictly
larger than one.  In this case, the events $\{ x_i \mbox{ is }
(m,l)-\mbox{regular} \}$, $i=1,2$, are independent.  We have
\begin{equation}
  \label{eq:rr.26}
{\cal A}^c \;\subset\; \bigcup_{x_1, \ldots, x_{R+1} \in
  \Lambda[x;L+\delta]}^{(l)} \bigcap_{i=1}^{R+1}
\Bigl\{ x_i \mbox{ is } (m,l)-\mbox{singular} \Bigr\},
\end{equation}
where the index $(l)$ reminds us that the union is over
collections of $R+1$ sites which are $l$-nonoverlapping.  Hence
\begin{eqnarray}
\PP\{ {\cal A}^c \} &\leq& \sum_{x_1, \ldots, x_{R+1} \in
  \Lambda[x; L+\delta]}^{(l)} \Bigl[
\PP\Bigl\{ x \mbox{ is } (m,l)-\mbox{singular} \Bigr\}
\Bigr]^{R+1} \nonumber \\
&\leq& \frac{[2(L+\delta)]^{d(R+1)}}{l^{p(R+1)}}
\; \leq\; \frac{1}{L^p} \frac{4^{d(R+1)}}{l^{(p - \alpha d)}}
\;\leq\; \frac{1}{2L^p}\;,
\end{eqnarray}
for $l$ large enough.   $\Box$

\begin{lema}
\label{varsovia2}
Let $l>\delta$, $\alpha >1$ and $p> 0$. Put $L = l^{\alpha}$. For
$\kappa$
and $b$ such that $0< \kappa < b/d$ take a collection of sites $\{x_i : i
=1,2,  \ldots ,R\}$ and introduce the event ${\cal B}_{\Delta}$ as in
  (\ref{eq:r20}) for $\Delta=e^{-l^\eta}$, with
  $0<\eta<b-\kappa d$. If $\aleph<\infty$ for some $a>\alpha(p+d)/\eta$,
then
\begin{equation}
\label{eq:r22}
\P \{  {\cal B}_\Delta \} \;\ge\; 1-{\frac{1}{2 L^p}}
\end{equation}
for $l$ larger than a certain threshold that depends on $\alpha$, $p$,
$\kappa$, $b$, $d$, $\eta$, $\que$.
\end{lema}

\proof We start with the following observations:
\begin{enumerate}
\item[(i)] For any collection of sites $\{x_i : i =1,2, \ldots ,R\}$
  \begin{equation}
    \label{eq:rr.28}
\ani_{\cup_{i=1}^{R} \Lambda[x_i,l^{\kappa}]} \;\subset\;
\bigcup_{i=1}^{R} \ani_{\Lambda[x_i;l^{\kappa} + \delta]}\;.
  \end{equation}
\item[(ii)] As a consequence
of (i),
\begin{equation}
\label{mas1}
{\cal B}_{\Delta}^{c} \;\subseteq\; \bigcup_{y \in
  \Lambda[x;L+\delta]} \, \Bigl\{
\sum_{\gamma \in \ani_{\Lambda[y;l^{\kappa}+\delta]}} \ln
\Bigl[\frac{1}{K_{\Delta}^{\jota}(\gamma)} \Bigr] \;>\;
\frac{l^b}{R} \Bigr\}.
\end{equation}
\item[(iii)] For any $\Delta>0$,
\begin{equation}
\label{mas2}
K_{\Delta}^{\jota}(\gamma) \;\geq\; (1 - e^{-\Delta})
\,e^{-\Delta w^{\jota}(\gamma)} \;>\; 0.
\end{equation}
\end{enumerate}

>From (\ref{mas2}) and the inequality
$\ln(1-e^{-\Delta}) \geq \ln (\Delta /2)$, valid for $\Delta$ small
enough, we obtain
\begin{eqnarray}
\label{mas3}
\sum_{\gamma \in \ani_{\Lambda[y;l^{\kappa}+\delta]}}\, \ln
\Bigl[\frac{1}{K_{\Delta}^{\jota}(\gamma)}
\Bigr] &\leq& \sum_{\gamma \in \ani_{\Lambda[y;l^{\kappa}+\delta]}} \, \Bigl[
\Delta  w^{\jota}(\gamma) - \ln(1 - e^{-
  \Delta})  \Bigr] \nonumber \\
&\leq& \sum_{z \in \Lambda[y;l^{\kappa}+\delta]} \, \sum_{ \gamma \ni z }
\,\Bigl[
\Delta  w^{\jota}(\gamma) - \ln(1 - e^{-
  \Delta})  \Bigr] \nonumber \\
&\leq& 4^d l^{\kappa d} \Bigl[ \Delta \sup_x \sum_{\gamma \ni x}
w^{\jota}(\gamma) - \Theta_{\ani} \ln \bigl( \frac{\Delta}{2} \bigr) \Bigr].
\end{eqnarray}
with $\Theta_{\ani} := \sup_x |\{\gamma \ni
x\}|  < \infty$.
Therefore, for $\Delta=e^{-l^{\eta}}$, $0< \eta < b - \kappa d$, and
assuming $l$ large enough in every step we get
\begin{eqnarray}
\label{mas4}
\P \Bigr\{ \sum_{\gamma \in \ani_{\Lambda[y;l^{\kappa}+\delta]}}\, \ln
\Bigl[\frac{1}{K_{\Delta}^{\jota}(\gamma)}
\Bigr] \;>\; \frac{l^b}{R} \Bigr\} &\leq& \P
\Bigl\{ \Delta\, \sup_x \sum_{\gamma \ni x} w^{\jota}(\gamma) > \frac{l^{b -
    \kappa
    d}}{4^d\, R} + \Theta_{\ani} \ln \bigl( \frac{\Delta}{2}
\bigr) \Bigr\}   \nonumber \\
&\leq& \P
\Bigl\{ \sup_x \sum_{\gamma \ni x} w^{\jota}(\gamma) > \frac{l^{b -
    \kappa
    d}}{4^{d+1}\,R\,\Delta} \Bigr\} \nonumber \\
&\leq& \P
\Bigl\{\ln \Bigl( 1 + \sup_x \sum_{\gamma \ni x} w^{\jota}(\gamma) \Bigr) >
\frac{l^{\eta}}{2} \Bigr\} \nonumber \\
&\leq& \frac{2^a}{l^{a \eta}} \,
\E \biggl[ \ln^{\textstyle a} \Bigl(1 + \sup_x \sum_{\gamma \ni x}
w^\jota(\gamma) \Bigr) \biggr],
\end{eqnarray}
where in the last step we have used Chebyshev inequality.  Finally,
\reff{mas1}, \reff{mas4} and the fact that $a \eta > \alpha (p+d)$
imply that
\begin{equation}
\label{mas5}
\P \{ {\cal B}_{\Delta}^c \} \;\leq\;
\frac{1}{L^p}\, \frac{2^{(a + 2d)} \aleph}{l^{a \eta - \alpha (p+d)}}
\;\leq\;  \frac{1}{2L^p},
\end{equation}
for $l$ large enough.  $\Box$

\subsubsection{Discussion on optimal choices}

We pause to discuss the possible choices for the parameters involved in
Lemmas \ref{geoest}, \ref{varsovia1} and \ref{varsovia2}.  The
parameter $a$ in Lemma \ref{varsovia2} involves explicitly the birth
rates of the animals.  We observe that the smaller the value of $a$
the weaker the condition on the birth rates.  Thus, the goal is to
choose the smallest possible value of this parameter compatible with
the hypotheses.  Let us first list the requirements on the different
parameters:
\begin{eqnarray}
\alpha &>& 1\;;\label{cond-1}\\
0\;\,<\;\, \nu &<& 1\;;\label{cond0}\\
\label{condicion1}
\max\{1, \nu +\theta_0 \} &<& \kappa \;\,<\;\, \frac{b}{d}
\;\,<\;\, \frac{\alpha \nu}{d}\;;\\
0\;\,<\;\; \eta &<& b - \kappa d\;;\label{cond1.6}\\
p &>& \alpha d\;;\label{cond1.7}\\
a &>& \frac{\alpha(p+d)}{\eta} \label{cond1.9}
\end{eqnarray}
Inequalities \reff{condicion1} tell us that for choice of $\kappa$ and
$b$ to be possible we must have that $\alpha\nu>d$ and hence that
[see \reff{eq:rr.30}]
\begin{equation}
\label{condicion3}
\theta_0 = \alpha(1-\nu) \;.
\end{equation}
This identity, in combination with the extreme and leftmost
inequalities in \reff{condicion1}, respectively implies that
\begin{equation}
\label{condicion4}
 \nu \;>\; \frac{\alpha d}{\alpha -d + \alpha d}\;,
\end{equation}
and
\begin{equation}
  \label{eq:rr.31}
  \kappa \;\ge\; \nu + \alpha(1-\nu)\;.
\end{equation}
Note that \reff{condicion4} yields, in view of \reff{cond0}, that
\begin{equation}
\label{condicion5}
\alpha \;>\; d\;.
\end{equation}
On the other hand, if we combine the rightmost inequality in
\reff{condicion1} with \reff{eq:rr.31} and \reff{cond1.6} we obtain
\begin{equation}
\label{condicion6}
\eta \;<\; \nu [ \alpha - d + \alpha d] - \alpha d \;<\; \alpha -d
\end{equation}
[the last inequality is due to \reff{cond0}].  Finally, we use this
last bound, in combination with \reff{cond1.7} and \reff{cond1.9} to
obtain a lower bound for the parameter $a$:
\begin{equation}
  \label{eq:rr.32}
 a \;>\; \frac{\alpha d (\alpha +1)}{\alpha - d} \qquad
\mbox{ with } \alpha > d\;.
\end{equation}

Small values of $a$ compatible with this restriction
are obtained when
\begin{equation}
\label{condicion7}
\alpha \;=\; d + \sqrt{d^2 + d}
\end{equation}
and satisfy the inequality
\begin{equation}
\label{condicion8}
a \;>\; 2 d^2 \Bigl( 1 + \sqrt{1+ {1 \over d}} + {1\over
  2d}\Bigr)\;.
\end{equation}
Reversing the preceding analysis, we see that the best strategy is to
take $\alpha$ and $a$ acording to \reff{condicion7} and
\reff{condicion8} and to choose $\nu$ and $p$ satisfying \reff{pf3}
and \reff{pf4}.  This yields the smallest possible choices of $a$ that
allow choices of $\eta$, $b$ and $\kappa$ respecting
the constraints \reff{cond-1}--\reff{cond1.9}.  This explains the
definition adopted for the parameters in the statement of the
change-of-scale Theorem \ref{porfin}.

\medskip

This is also a good opportunity to comment on choices for time
scales.  They are determined by the demands impossed throughout the
proofs of Sublemma \ref{sl.second} and Lemma \ref{varsovia2}.  There
are two major constraints:  First, the scales $\widehat\delta(l)$ and
$\widehat\tau(l)$ must be truly intermediate in the sense that
\begin{equation}
\label{con1}
l \;\prec\; \widehat\delta(l) \;\prec\; L=l^\alpha
 \qquad\hbox{and} \qquad
T(l) \;\prec\; \widehat{\tau}(l) \;\prec\; T(L)=T(l^\alpha)
\;,
\end{equation}
where $f(l)\prec g(l)$ means that $f(l)/g(l)$ tends to zero as $l$
tends to infinity.  Second, the exponent $h(l)$ defined in
\reff{qudefe}/\reff{haga} must grow with $l$.  This implies that
\begin{equation}
  \label{eq:rr.20.0}
  \widehat\tau(l) \;\prec\; e^{\widehat\delta(l)}\;.
  \end{equation}
  The combination of the two preceding displays yields the relations
\begin{equation}
  \label{eq:rr.20.1}
  T(l) \;\prec\; \widehat\tau(l) \;\prec\; e^{\widehat\delta(l)}
  \;\prec\; e^{l^\alpha}
\end{equation}
which shows that at best the time scale can grow as a stretched
exponential.

\subsubsection{Conclusion of the proof of Theorem \protect\ref{porfin}}
We prove \reff{pipio} by induction in $k$.  The case $k=0$ follows
from hypothesis \reff{pf5} because $(m_0,L)$-regularity gives
$(m_{\infty},L)$-regularity for $m_{\infty}<m_0$.  Suppose now
 \reff{pipio} is true for some $k>0$.  Lemmas \ref{varsovia1} and
\ref{varsovia2} ---with the replacements $(m,l) \rightarrow (m_k,L_k)$
and $L \rightarrow L_{k+1}$--- and the preceding discussion on optimal
choices imply that if
\begin{equation}
  \label{eq:rr.33}
\P\Bigl\{ x \mbox{ is } (m_k,L_k)-\mbox{regular} \Bigr\} \;\geq\; 1 -
\frac{1}{L_k^p}\;,
\end{equation}
then
\begin{equation}
\label{pastis1}
\P \{{\cal A} \cap {\cal B}_{\Delta} \} \;\geq\; 1 -
\frac{1}{L_{k+1}^p}
\end{equation}
for $l$, and hence $L_k$, sufficiently large.
On the other hand, Lemma \ref{geoest} ---with
$(m,l) \rightarrow (m_k,L_k)$ and $(M,L) \rightarrow (m_{k+1},L_{k+1})$---
implies that if $L_k^{-\theta} < m_k < m_0$ for $0< \theta < \theta_0
= \alpha(1-\nu)$, then
\begin{equation}
\label{pastis2}
{\cal A} \cap {\cal B}_{\Delta}  \;\subseteq\;
\Bigl\{  x \mbox{ is } (m_{k+1},L_{k+1})-\mbox{regular} \Bigr\}\;,
\end{equation}
with
\begin{equation}
  \label{eq:rr.34}
m_{k+1} \;\geq\; m_k - \frac{a_0}{L_k^{\theta_0}} \;\geq\;
\frac{1}{L_{k+1}^{\theta}}
\end{equation}
for $l$, and hence $L_k$, sufficiently large.
>From (\ref{pastis1}) and (\ref{pastis2}) we get that
\begin{equation}
\label{pastis3}
\P\Bigl\{ x \mbox{ is } (m_k,L_k)-\mbox{regular} \Bigr\}
\;\geq\; 1 - \frac{1}{L_k^p}
\end{equation}
for all natural $k$, for $l$ suffciently large.

To conclude we must check that
\begin{equation}
  \label{eq:rr.35}
m_{\infty} \;\leq\; m_k \quad , \quad \mbox{for } k=1,2, \ldots
\end{equation}
>From \reff{eq:rr.34} we have
\begin{equation}
  \label{eq:rr.36}
m_k \;\geq\; m_0 - a_0 \sum_{j=0}^{k-1} \, \frac{1}{L_j^{\theta_0}}
\;=\;  m_0 - a_0
\sum_{j=0}^{k-1} \, \Bigl(\frac{1}{L_0^{\theta_0}}\Bigr)^{\alpha^j}\;.
\end{equation}
Therefore, \reff{eq:rr.35} is verified if $L_0$ is chosen so large
that
\begin{equation}
  \label{eq:rr.37}
a_0 \sum_{j=0}^{\infty} \, \Bigl(\frac{1}{L_0^{\theta_0}}
\Bigr)^{\alpha^j} \;<\; m_0 - m_{\infty}\;. \quad \Box
\end{equation}

\subsection{A last probabilistic estimate:  The choice of initial
  scale}\label{s.fin} To finish the proof of the key Lemma \ref{main}
we show that for every real $\widetilde L$ (though we are interested
in $\widetilde L$ large) there exists
$L_0\ge \widetilde L$ such that \reff{pf5} holds.
It is only here
that hypothesis \reff{eq:r3} is needed.  Let us suppose that
a function $S:\ani\to [1,\infty)$ has been selected and consider the
resulting quantity $\Psi^\jota$ defined in \reff{eq:r30}.

\begin{lema} \label{lem:rr.1}
  Given $\widetilde L, p>0$, for each $\rho\in[0,1)$ there exist
  $m(\rho), \varepsilon(\rho)>0$ and $L(\rho)\ge\widetilde L$ such
  that
\begin{equation}
\label{pf5bis}
\P\Bigl\{ \Psi^\jota > \rho \Bigr\}\;<\; \varepsilon(\rho)
\quad\Longrightarrow\quad
\P\Bigl\{ x \mbox{ is } (m(\rho) ,L(\rho))-\mbox{regular} \Bigr\}
\;\geq\; 1 - \frac{1}{L(\rho)^p}\;.
\end{equation}
\end{lema}

\proof
By Theorem 5.1 in Fern\'{a}ndez, Ferrari and Garcia
(2001)\nocite{ffg01}, the condition $\Psi^\jota\le \rho$ implies that
there exist $c(\rho), \widetilde m(\rho)>0$ such that
\begin{equation}
  \label{eq:rr.40}
  G_{B_{L + \delta,T}(X)}^{\jota}(X,Y) \;\le\; c(\rho)\,
\exp\Bigl\{-\widetilde m(\rho)\,
\max(\norm{x-y}\,, -t_Y) \Bigr\}\;,
\end{equation}
for all $Y=(y,t_Y)$ with $t_Y\le 0$.  Taking, for instance,
$m(\rho)=\widetilde m(\rho)/2$ we conclude that there exists
$\widetilde L(\rho)$ such that
\begin{equation}
  \label{eq:rr.41}
  \Bigl\{ x \mbox{ is } (m(\rho) ,L)-\mbox{regular} \Bigr\}
\;\supset\; \Bigl\{ \Psi^\jota \le \rho \Bigr\}
\end{equation}
for all $L\ge\widetilde L(\rho)$.  Implication \reff{pf5bis} follows by
choosing $L(\rho)=\max [\widetilde L\,,\, \widetilde L(\rho)]$ and
$\varepsilon(\rho)=L(\rho)^{-p}$.  $\square$
\medskip

Chebyshev inequality shows that \reff{eq:r3} is a sufficient condition
for the validity of the left-hand side of \reff{pf5bis}, as long as
$\varepsilon\le \rho\,\varepsilon(\rho)$, for any given choice of
$\rho$.  This observation, the preceding Lemma \ref{lem:rr.1} and
Theorems \ref{porfin} and \ref{behavior} constitute the proof of the
main Lemma \ref{main}.  We observe that, in our proof, $q=
\frac{1}{\nu}$ and $q_0 = \frac{a[\alpha - d + \alpha d]}{\alpha
  d[\alpha + a +1]}$ with $\alpha = d + \sqrt{d^2 + d}$.

\section*{Acknowledgements}

The essential part of this work was done while G.R.G.\ was a visiting
scholar at the Instituto de Matem\'{a}tica e Estat\'{\i}stica,
Universidade de S\~ao Paulo.  He would like to acknowledge his warm
gratitude to this institution and to FAPESP, the funding agency.  He
also thanks the Laboratoire de Math\'{e}matiques Rapah\"el Salem, UMR
6085, CNRS Univer\-si\-t\'{e} de Rouen for inviting him during the
completion of his work.  R.F.\ thanks the aforementioned Instituto de
Matem\'{a}tica e Estat\'{\i}stica de la Universidade de S\~ao Paulo and
the Newton Institute for the Mathematical Sciences at Cambridge
University, for hospitality during his work in this paper.  The
authors wish to thank J.~van den~Berg, G.~Grimmett, F.~den Hollander
and A.~Klein for valuable discussions.  This work received the support
of FAPESP, CNPq, FINEP (N\'{u}cleo de Excel\^encia \emph{Fen\^omenos
  cr\'{\i}ticos em probabilidade e processos estoc\'{a}sticos},
PRONEX-177/96) and an agreement USP-COFECUB (Project UC 75/01:
\emph{Les Processus Al\'{e}atoires et la M\'{e}canique Statistique
  Math\'{e}matique}).


\begin{thebibliography}{10}

\bibitem{bkl96}
A.~J. Baddeley, W.~S. Kendall, and M.~N.~M. Van~Lieshout.
\newblock Quermass-interaction porcesses, 1996.
\newblock Preprint.

\bibitem{badlie95}
A.~J. Baddeley and M.~N.~M. van Lieshout.
\newblock Area-interaction point processes.
\newblock {\em Ann. Inst. Statist. Math.}, 47(4):601--619, 1995.

\bibitem{bezgri91}
C.~Bezuidenhout and G.~Grimmett.
\newblock Exponential decay for subcritical contact and percolation processes.
\newblock {\em Ann. Probab.}, 19(3):984--1009, 1991.

\bibitem{brikup96}
J.~Bricmont and A.~Kupiainen.
\newblock High temperature expansions and dynamical systems.
\newblock {\em Comm. Math. Phys.}, 178(3):703--732, 1996.

\bibitem{camkleper91}
M.~Campanino and A.~Klein.
\newblock Decay of two-point functions for {$(d+1)$}-dimensional percolation,
  {I}sing and {P}otts models with {$d$}-dimensional disorder.
\newblock {\em Comm. Math. Phys.}, 135(3):483--497, 1991.

\bibitem{camkle91}
M.~Campanino, A.~Klein, and J.~F. Perez.
\newblock Localization in the ground state of the {I}sing model with a random
  transverse field.
\newblock {\em Comm. Math. Phys.}, 135(3):499--515, 1991.

\bibitem{dob96}
R.~L. Dobrushin.
\newblock Perturbation methods of the theory of {G}ibbsian fields.
\newblock In {\em Lectures on probability theory and statistics (Saint-Flour,
  1994)}, pages 1--66. Springer, Berlin, 1996.

\bibitem{ffg98}
R.~Fern{\'{a}}ndez, P.~A. Ferrari, and N.~L. Garcia.
\newblock Measures on contour, polymer or animal models. {A} probabilistic
  approach.
\newblock {\em Markov Process. Related Fields}, 4(4):479--497, 1998.
\newblock I Brazilian School in Probability (Rio de Janeiro, 1997).

\bibitem{ffg01}
R.~Fern{\'{a}}ndez, P.~A. Ferrari, and N.~L. Garcia.
\newblock Loss network representation of {P}eierls contours.
\newblock {\em Ann. Probab.}, 29:902--937, 2001.

\bibitem{ffg00}
R.~Fern{\'{a}}ndez, P.~A. Ferrari, and N.~L. Garcia.
\newblock Perfect simulation for interacting point processes, loss networks and
  ising models.
\newblock {\em Stoch.\ Proc. Appl.}, 102:63--88, 2002.

\bibitem{fg98}
P.~A. Ferrari and N.~L. Garcia.
\newblock One-dimensional loss networks and conditioned ${M}/{G}/\infty$
  queues.
\newblock {\em J. Appl. Probab.}, 35(4):963--975, 1998.

\bibitem{giemae96}
G.~Gielis and C.~Maes.
\newblock Percolation techniques in disordered spin flip dynamics: relaxation
  to the unique invariant measure.
\newblock {\em Comm. Math. Phys.}, 177(1):83--101, 1996.

\bibitem{gri95}
G.~Grimmett.
\newblock The stochastic random-cluster process and the uniqueness of
  random-cluster measures.
\newblock {\em Ann.\ Probab.}, pages 1461--1450, 1995.

\bibitem{kel91}
F.~P. Kelly.
\newblock Loss networks.
\newblock {\em Ann. Appl. Probab.}, 1(3):319--378, 1991.

\bibitem{ken98}
W.~S. Kendall.
\newblock Perfect simulation for the area-interaction point process.
\newblock In L.~Accardi and C.~C. Heyde, editors, {\em Probability Towards
  2000}, pages 218--234. Springer, 1998.

\bibitem{kle94}
A.~Klein.
\newblock Extinction of contact and percolation processes in a random
  environment.
\newblock {\em Ann. Probab.}, 22, 1994.

\bibitem{kle95}
A.~Klein.
\newblock Multiscale analysis in disordered systems: {P}ercolation and contact
  process in a random environment.
\newblock Preprint, 1995.

\bibitem{meeroy96}
R.~Meester and R.~Roy.
\newblock {\em Continuum percolation}, volume 119 of {\em Cambridge Tracts in
  Mathematics}.
\newblock Cambridge University Press, Cambridge, 1996.

\bibitem{mol01}
J.~M{\o}ller.
\newblock A review of perfect simulation in stochastic geometry.
\newblock In {\em Selected Proceedings of the Symposium on Inference for
  Stochastic Processes}, pages 333--55. IMS Lecture Notes \& Monographs Series,
  Volume 37, 2001.

\bibitem{molwaa04}
J.~M{\o}ller and R.~P. Waagepetersen.
\newblock {\em Statistical inference and simulation for spatial point
  processes}, volume 100 of {\em Monographs on Statistics and Applied
  Probability}.
\newblock Chapman \& Hall/CRC, Boca Raton, FL, 2004.

\bibitem{str75}
D.~J. Strauss.
\newblock A model for clustering.
\newblock {\em Biometrika}, 62(2):467--475, 1975.

\bibitem{ber96}
J.~van~den Berg.
\newblock A note on disjoint-occurrence inequalities for marked {P}oisson point
  processes.
\newblock {\em J. Appl. Probab.}, 33(2):420--426, 1996.

\bibitem{dre87}
H.~von Dreifus.
\newblock On the effects of randomness in ferromagnetic models and
  {Schr\"odinger} operators.
\newblock Phd. Thesis, New York University, 1987.

\bibitem{widrow70}
B.~Widow and J.~S. Rowlinson.
\newblock New model for the study of liquid-vapor phase transitions.
\newblock {\em J. Chem. Phys.}, 52:1670--1684, 1970.

\end{thebibliography}

\end{document}